\documentstyle[12pt]{article}

\author{H.Falomir$^1$
\and
R.E.Gamboa Sarav\'\i $^1$
\and
M.A.Muschietti$^2$
\and
E.M.Santangelo$^1$
\and
J.E.Solomin$^2$ \\ \hfill
\and
{\Large Facultad de Ciencias Exactas, UNLP}\\
$^1$Departamento de F\'{\i}sica\\
$^2$Departamento de Matem\'atica}
\title{Determinants of elliptic boundary problems for Dirac operators
I.
Local boundary conditions.
\thanks{Partially supported by CONICET, Argentina. } }
\date{ }

\def\dfrac#1#2{{\displaystyle {#1 \over #2}}}
\def\binom#1#2{{#1 \choose #2}}
\def\QATOP#1#2{{#1 \atop #2}}
\def\QDATOP#1#2{{\displaystyle {#1 \atop #2}}}

\catcode`\@=11
\def\qed{\ifhmode\unskip\nobreak\fi\ifmmode\ifinner\else
\hskip5\p@\fi\fi
\hbox{\hskip5\p@\vrule width4\p@ height6\p@ depth1.5
\p@\hskip\p@}}

\begin{document}

\maketitle
\begin{abstract}
We study functional determinants for Dirac operators on manifolds
with
boundary and discuss the ellipticity of boundary problems by using
the
Calder\'on projector. We give, for local boundary conditions, an
explicit
formula relating these determinants to the corresponding Green
functions. We
finally apply this result to the case of a bidimensional disk under
bag-like
conditions.
\end{abstract}
\eject

\section{Introduction}

It is well known that functional determinants have wide application
in
Quantum and Statistical Physics. Typically, one faces the necessity
of
defining a regularized determinant for elliptic differential
operators. In
this context, the Dirac first order differential operator plays a
central role.

Seeley's construction of complex powers of elliptic differential
operators
provides a powerful tool to regularize such determinants: the so
called \linebreak
$\zeta$-function method \cite{Haw}.

In the case of boundaryless manifolds, this construction has been
largely
studied and applied (see, for instance, \cite{Annals} and references
therein).

For manifolds with boundary, the study of complex powers was
performed in  \cite{Seeley1,Seeley2} for the case of local boundary
conditions,
while for the case of nonlocal conditions, this task is still in
progress
(see, for example, \cite{SG}.)

In general, the regularized determinant turns out to be nonlocal and
so, it
cannot be expressed in terms of just a finite number of Seeley's
coefficients. However, such determinant can always be obtained from
the
Green function in a finite number of steps involving these
coefficients. For
boundaryless manifolds this was proved in \cite{Jour.1}, while for a
particular type of local boundary conditions the procedure was
introduced in
\cite{Jour.2}.

The aim of this paper is to give a rigorous proof of the validity of
this
assertion in the case of the Dirac operator under general local
elliptic
boundary conditions. In so doing, the explicit relationship between
determinants and the corresponding Green functions will  be derived.

Dirac operators defined on manifolds with boundaries have been the
subject
of a vast literature (see, for instance, \cite{NT,Moreno} and
references therein), mainly concerning anomalies and index theorems.
In
these papers, the emphasis was put on nonlocal boundary conditions of
the
type introduced in \cite{APS}. We leave for a forthcoming publication
the
treatment of such conditions.

Less attention was devoted to local boundary conditions in physical
literature. In fact, some problems do not even admit such conditions
owing
to {\it topological obstructions }. We find enlightening to make a
detailed
discussion of this point since, to our knowledge, it has not been
done in
this context.

The outline of this paper is as follows:

In Section 2 we study elliptic boundary problems for Dirac operators,
by
means of the Calder\'on projector \cite{LibroAmarillo}. We discuss in
detail
the topological obstructions, and consider a class of local boundary
conditions giving rise to elliptic boundary systems.

In Section 3 we review, for the sake of completeness, Seeley's
construction
for complex powers of pseudo - differential operators on manifolds
with
boundary.

In Section 4 we present the main result of the paper: A formula
relating the
determinant of the Dirac operator with its Green function is
established.

In Section 5 an explicit computation of the determinant of a Dirac
operator
in a bidimensional disk with bag-like boundary conditions is given.

\section{The Calder\'on projector and Elliptic boundary problems}

Throughout this paper we will be concerned with boundary value
problems
associated to first order elliptic operators
\begin{equation}
\label{OP}D:C^\infty (M,E)\rightarrow C^\infty (M,F),
\end{equation}
where $M$ is a compact connected sub manifold of ${\bf R}^\nu $ of
codimension zero with smooth boundary $\partial M$, and E and $F$ are
$k$%
-dimensional complex vector bundles over $M.$

In a collar neighborhood of $\partial M$ in $M,$ we will take
coordinates $%
\bar x=(x,t)$, with $t$ the inward normal coordinate and $x$ local
coordinates for $\partial M$ (that is, $t>0$ for points in $M$
$\setminus $$%
\partial M$ and $t=0$ on $\partial M$ ), and conjugated variables
$\bar \xi
=(\xi ,\tau )$.

As stated in the Introduction, we will mainly consider  the Euclidean
Dirac operator. Let us
recall that the free Euclidean Dirac operator $\not \! \partial $ is
defined as
\begin{equation}
i\not \! \partial =\sum_{\mu =0}^{\nu -1}i\gamma _\mu \frac \partial
{\partial
x_\mu },
\end{equation}
where the matrices $\gamma _\mu $ satisfy
\begin{equation}
\gamma _\mu \gamma _\alpha +\gamma _\alpha \gamma _\mu =2\delta _{\mu
\alpha
},
\end{equation}
and that, given a gauge potential $A=\{A_\mu ,\ \mu =0,...,\nu -1\}$
on $M$,
the coupled Dirac operator is defined as
\begin{equation}
\label{cul}D(A)=i\not \! \partial +\not \! \! A
\end{equation}
with $\not \!\! A=\sum\limits_{\mu =0}^{\nu -1}\gamma _\mu A_\mu .$

In Section 5, explicit computations are carried out for $\nu =2.$ We
take, for
this case, the representation for the $\gamma ^{\prime }s$ given by
the
Pauli matrices \linebreak \{$\sigma _j,\ j=1,2,3\}$

\begin{equation}
\label{gm}
\begin{array}{c}
\gamma _0=\sigma _1=\left(
\begin{array}{cc}
0 & 1 \\
1 & 0
\end{array}
\right) ~,\qquad \gamma _1=\sigma _2=\left(
\begin{array}{cc}
0 & -i \\
i & 0
\end{array}
\right) ~, \\
\\
\qquad \gamma _5=-i\gamma _0\gamma _1=\sigma _3=\left(
\begin{array}{cc}
1 & 0 \\
0 & -1
\end{array}
\right) .
\end{array}
\end{equation}
In this representation, fermions with positive chirality are of the
form $\binom{%
\bullet }0$, and those with negative chirality of the form $\binom
0{\bullet
}. $ So, the free Dirac operator can be written as:
\begin{equation}
\label{dl}i\not \! \partial =i\ (\gamma _0\ \partial _0+\gamma _1\
\partial
_1)=2i\ \left(
\begin{array}{cc}
0 & \frac \partial {\partial z} \\
\frac \partial {\partial z^{*}} & 0
\end{array}
\right) ,
\end{equation}
where $z=x_0+i\ x_1.$

One of the most suitable tools for studying boundary problems is the
Calder\'on projector $Q$ \cite{LibroAmarillo}. For the case we are
interested in, $D$ of order 1 as in (\ref{OP}), $Q$ is a (not
necessarily
orthogonal) projection from [$L^2(\partial M,E_{/\partial M})]$ onto
the
subspace $\{(T\varphi \ /\varphi \in \ker (D)\},$ being $T:C^\infty
(M,E)\rightarrow C^\infty (\partial M,E_{/\partial M})$ the trace
map.

Given any fundamental solution $K(\bar x,\bar y)$ of $D,$ the
projector $Q$
can be constructed in the following way: for $f\in C^\infty (\partial
M,E_{/\partial M})$ , one gets $\ \varphi \in \ker (D)$ by means of a
Green
formula involving $K(\bar x,\bar y)$, and takes the limit of $\varphi
$ for $\bar x\rightarrow \partial M$.

As shown in \cite{LibroAmarillo}, $Q$ is a zero-th order pseudo
differential
operator and its principal symbol $q(x;\xi ),$ that depends only on
the
principal symbol of $D,\sigma _1(D)=a_1(x,t;\xi ,\tau ),$ turns out
to be
the $k\times k$ matrix
\begin{equation}
\label{q1}q(x;\xi )=\frac 1{2\pi i}\int_\Gamma \left(
a_1^{-1}(x,0;0,1)\
a_1(x,0;\xi ,0)-z\right) ^{-1}dz,
\end{equation}
where $\Gamma $ is any simple closed contour oriented clockwise and
enclosing all poles of the integrand in $Im(z)<0.$

Although Q is not uniquely defined, since one can take any
fundamental solution $K$ of $D$ to construct it, the image of $Q$ and
its
principal symbol $q(x;\xi )$ are independent of the choice of $K$
\cite
{LibroAmarillo}.

We find enlightening to compute the principal symbol of the
Calder\'on
projector for the Dirac operator as in (\ref{cul}) directly from the
definition of $Q$, instead of using (\ref{q1}).

Let $K(\bar x,\bar y)$ be a fundamental solution of the Dirac
operator $D(A)$
in a neighborhood of the region M$\subset {\bf R}^\nu $, i.e.

\begin{equation}
D^{\dagger }(A)K^{\dagger }(\bar x,\bar y)=\delta (\bar x-\bar y).
\end{equation}
We can write
\begin{equation}
\label{po}K(\bar x,\bar y)=K_0(\bar x,\bar y)+R(\bar x,\bar y)
\end{equation}
where $K_0(\bar x,\bar y)$ is a fundamental solution of
$i\!\!\!\!\!\!\not\!\!\!\!\!\!\partial$
and $\vert R(\bar x,\bar y)\vert $ is \linebreak  O(1/$\vert \bar
x-\bar y\vert ^{\nu -2})$
for $\bar x-\bar y\sim 0.$ Since $(i\! \! \! \not \! \! \partial
)^2=-\bigtriangleup $,
it is easy to obtain $K_0(\bar x,\bar y)$ from the well known
fundamental
solution of the Laplacian, so
\begin{equation}
\label{pot}K_0(\bar x,\bar y)=-\ i\ \frac{\Gamma (\nu /2)}{2\ \pi
^{\nu /2}}%
\ \ \frac{(\bar {\not\! x}-\bar {\not\! y})}{\vert \bar x-\bar y\vert
^\nu }.
\end{equation}
For $f$ a smooth function on $\partial M,$%
\begin{equation}
\label{poto}Q f(x)=-i\lim \limits_{\bar x\rightarrow \partial
M}\int_{\partial M}K(\bar x,y)\ \not \! n\ f(y)\ d\sigma _y,
\end{equation}
where $\not \!\! n=\sum_l\gamma _l\ n_l,$ and $n=(n_l)$ is the
unitary outward
normal vector on $\partial M.$ Note that, if $f=T\varphi $, with
$\varphi
\in \ker (D)$, the Green formula yields $Qf=f$, as required.

{}From (\ref{po}), (\ref{pot}) and (\ref{poto}) one gets
\begin{equation}
\label{101}Q f(x)=\frac 12f(x)-i\ P.V.\int_{\partial M}K(x,y)\ \not
\! n\
f(y)\ d\sigma _y
\end{equation}
In order to see that the principal value in (\ref{101}) makes sense,
note
that
\begin{equation}
\label{102}P.V.\int_{\partial M}K_0(x,y)\ \not \! n\ f(y)\ d\sigma _y
\end{equation}
is meaningful since
\begin{equation}
\label{103}-i\int_{\partial M}K_0(\bar x,y)\ \not \! n\ \ d\sigma
_y=Id
_{k\times k},\ \forall \bar x\in M.
\end{equation}
In fact, any $\varphi (\bar x)\equiv constant$ is a solution of
$i\not \! \partial \varphi =0$ in $M.$

Then
\begin{equation}
\label{104}
\begin{array}{c}
Q f(x)=\frac 12f(x)-i\ P.V.\int_{\partial M}K_0(x,y)
\not \! n\ f(y)\ d\sigma _y\  \\  \\
-i\ \int_{\partial M}R(x,y)\ \not \! n\ f(y)\ d\sigma _y
\end{array}
\end{equation}
For the calculus of the principal symbol, we write the second term in
the
r.h.s. of (\ref{104}) in local coordinates on $\partial M$ ,
\begin{equation}
\label{105}\ \ -iP.V.\int_{{\bf R}^{\nu -1}}\frac{\Gamma (\nu /2)}{2\
\pi
^{\nu /2}}\ \frac{(x-y)_j}{\vert x-y\vert ^\nu }\ \gamma _{j\ }\gamma
_n f(y)\
dy=\frac 1{2\ }\ \gamma _{j\ }\gamma _n\ {\bf R}_j{\bf (}f)(x),
\end{equation}
where ${\bf R}_j(f)$ is the $j$-$th$ Riesz transform of $f$. The
symbol of the
operator in (\ref{105}) is (see for example \cite{Stein})
\begin{equation}
\frac 12i\ \gamma _{j\ }\gamma _n\frac{\xi _j}{\vert\xi\vert}\ =\frac
12i\
\frac{\not \! \xi }{\vert \xi \vert }\ \not \! n.
\end{equation}
The last term in the r.h. side of (\ref{104}) is a pseudo -
differential
operator of order $\leq -1,$ because of the local behavior of
$R(x,y)$, and
then it does not contribute to the calculus of the principal symbol
we are
carrying out. Then, coming back to global coordinates, we finally
obtain

\begin{equation}
\label{q}q(x;\xi )=\frac 12(Id_{k\times k}+i\ \frac{\not \! \xi }
{\vert \xi \vert }\ \not \! n).
\end{equation}
In order to get the $rank$ of this matrix, note that
\begin{equation}
\label{rango}
\begin{array}{c}
q(x;\xi )\ q(x;\xi )=q(x;\xi ) \\
\\
tr\ q(x;\xi )={k/2},
\end{array}
\end{equation}
and consequently $rank\ q(x;\xi )=k/2.$

In particular for $\nu =2$ and the $\gamma ^{\prime }s$ as in
(\ref{gm}), we
obtain
\begin{equation}
\label{q2}q(x;\xi )=\left(
\begin{array}{cc}
H(\xi ) & 0 \\
0 & H(-\xi )
\end{array}
\right)
\end{equation}
$\forall x\in \partial M,$ with $H(\xi )$ the Heaviside function.

\bigskip\

According to Calder\'on \cite{LibroAmarillo}, elliptic boundary
conditions
can be defined in terms of $q(x;\xi )$, the principal symbol of the
projector $Q.$

{\bf Definition} 1{\sc :}

Let us assume that the $rank$ of $q(x;\xi )$ is a constant $r$ ( as
is
always the case for $\nu \geq 3$ \cite{LibroAmarillo}).

A zero order pseudo differential operator $B:[L^2(\partial
M,E_{/\partial
M})]\rightarrow $ $[L^2(\partial M,G)],$ with $G$ an $r$ dimensional
complex
vector bundle over $\partial M,$ gives rise to an {\it elliptic
boundary
condition }for a first order operator $D$ if, $\forall \xi :\vert \xi
\vert
\geq 1,$
\begin{equation}
\label{erre}rank(b(x;\xi )\ q(x;\xi ))=rank(q(x;\xi ))=r\ ,
\end{equation}
where $b(x;\xi )$ coincides with the principal symbol of $B$ for
$\vert \xi
\vert \geq 1.$

In this case we say that

\begin{equation}
\label{BoundaryProblem}\left\{
\begin{array}{c}
D\varphi =\chi\ \
\rm{ in }\ M \\  \\
BT\varphi =f\ \rm{ on }\ \partial M
\end{array}
\right.
\end{equation}
is an {\it elliptic boundary problem}, and denote by $D_B$ the
closure of $D$
acting on the sections $\varphi $ $\in C^\infty (M,E)$ satisfying $%
B(T\varphi )=0.$

In particular, condition (\ref{erre}) implies that, for each $s\in
{\bf R},$
the image of $B\circ Q$ is a closed subspace of the Sobolev space $%
H^s(\partial M,G)$ having finite codimension \cite{LibroAmarillo}.

An elliptic boundary problem as (\ref{BoundaryProblem}) has a
solution $%
\varphi \in H^1(M,E)$ for any $(\chi ,f)$ in a subspace of
$L^2(M,E)\times
H^{1/2}(\partial M,G)$ of finite codimension. Moreover, this solution
is
unique up to a finite dimensional kernel \cite{LibroAmarillo}. In
other
words, the operator
\begin{equation}
(D,BT):H^1(M,E)\rightarrow L^2(M,E)\times H^{1/2}(\partial M,G)
\end{equation}
is Fredholm.

For $\nu =2$, Definition 1 does not always apply. For instance, for
the two
dimensional chiral Euclidean Dirac operator
\begin{equation}
\label{cdo}D=2i\frac \partial {\partial z^{*}}\ ,
\end{equation}
acting on sections with positive chirality and taking values in the
subspace
of sections with negative one, it is easy to see from (\ref{q2}) that
\begin{equation}
\label{qb}q(x;\xi )=H(\xi ).
\end{equation}
Then, the $rank$ of $q(x;\xi )$ is not constant. In fact,
\begin{equation}
rank\ q(x;\xi )=\left\{ \QDATOP{0\quad \rm{ if }\ \xi <0}{1\quad \rm{
if }\ \xi >0}\right. .
\end{equation}
However, for the (full) two dimensional Euclidean Dirac operator
\begin{equation}
\label{fd2}D(A)=\left(
\begin{array}{cc}
0 & D^{\dagger } \\
D & 0
\end{array}
\right)
\end{equation}
we get from (\ref{rango}) that $rank\ q(x;\xi )=2/2=1$ $\forall \xi
\neq 0$, and so Definition 1 does apply.

When $B$ is a local operator, Definition 1 yields the classical local
elliptic boundary conditions, also called Lopatinsky-Shapiro
conditions (see
for instance \cite{Hor}) .

For Euclidean Dirac operators on ${\bf R}^\nu ,$ $E_{/\partial
M}=\partial
M\times {\bf C}^k,$ and local boundary conditions arise when the
action of $%
B $ is given by the multiplication by a $\frac k2\times k$ matrix of
functions defined on $\partial M.$

Owing to {\it topological obstructions, }chiral Dirac operators in
even
dimensions do not admit{\it \ }local elliptic boundary conditions
(see for
example \cite{Booss-B}). For instance, in four dimensions, by
choosing the $
\gamma $-matrices at \linebreak $x=(x_1,x_2,x_3)\in \partial M$ as
\begin{equation}
\gamma _n=i\left(
\begin{array}{cc}
0 & Id_{2\times 2} \\
-Id_{2\times 2} & 0
\end{array}
\right) \ \ \ \rm{and}\qquad \gamma _j=\left(
\begin{array}{cc}
0 & \sigma _j \\
\sigma _j & 0
\end{array}
\right) \ \rm{for }\ \ j=1,2,3,
\end{equation}
the principal symbol of the Calder\'on projector (\ref{q}) associated
to
the full Dirac
operator turns out to be
\begin{equation}
q(x;\xi )=\frac 12\left(
\begin{array}{cc}
Id_{2\times 2}+\dfrac{\xi .\sigma }{\vert \xi \vert } & 0 \\
0 & Id_{2\times 2}-\dfrac{\xi .\sigma }{\vert \xi \vert }
\end{array}
\right) .
\end{equation}
Thus, from the left upper block, one gets for the chiral Dirac
operator
\begin{equation}
q_{ch}(x;\xi )=\frac 12\left(
\begin{array}{cc}
1+\dfrac{\xi _3}{\vert \xi \vert } & \dfrac{\xi _1-i\xi _2}{\vert \xi
\vert } \\
  &  \\
\dfrac{\xi _1+i\xi _2}{\vert \xi \vert } & 1-\dfrac{\xi _3}{\vert \xi
\vert }
\end{array}
\right) .
\end{equation}
So $q_{ch}(x;\xi )$ is a hermitian idempotent 2$\times 2$ matrix with
$%
rank=1.$ If one had a local boundary condition with principal symbol
$%
b(x)=(\beta _1(x),\beta _2(x))$, according to Definition 1, it should
be $%
rank(b(x)\ q_{ch}(x;\xi ))=1,\ \forall \xi \neq 0.$ However, it is
easy to
see that for
\begin{equation}
\xi _1=\frac{-2\beta _1\beta _2}{\beta _1^2+\beta _2^2},\qquad \xi
_2=0\quad
\rm{and\quad }\xi _3=\frac{\beta _2^2-\beta _1^2}{\beta _1^2+\beta
_2^2},
\end{equation}
$rank(b(x)\ q_{ch}(x;\xi ))=0.$ Equivalently, this means that $\
q_{ch}(x;\xi )$ is not deformable through idempotents of $rank=1$ to
a $\xi
- $independent matrix function . This is an example of the above
discussed
topological obstructions.

Nevertheless, it is easy to see that local boundary conditions can be
defined for the full, either free or coupled, Euclidean Dirac
operator

$$
D(A)=\left(
\begin{array}{cc}
0 & D^{\dagger } \\
D & 0
\end{array}
\right)
$$
on $M.$ For instance, we see from (\ref{q2}) and (\ref{erre})
that for $\nu =2$, the
operator $B$ defined as
\begin{equation}
\label{betas}
B\left( \QATOP{f}{g}\right) =(\beta _1(x),\beta _2(x))\left(
\QATOP{f}{g}%
\right)
\end{equation}
yields a local elliptic boundary condition for every couple of
nowhere
vanishing functions $\beta _1$$(x)$ and $\beta _2(x)$ on $\partial
M.$

A type of non-local boundary conditions, to be consider in a
forthcoming
publication, is related to the ones defined and analyzed by M. Atiyah
,
V. Patodi and I. Singer in \cite{APS} for a wide class of first order
Dirac-like operators, including the Euclidean chiral case. Near
$\partial M$
such operators can be written as
\begin{equation}
\sigma \ (\partial _t+A),
\end{equation}
where $\sigma $ is an isometric bundle isomorphism $E\rightarrow F,$
and $A:L^2(\partial M,E_{/\partial M})\rightarrow $ $L^2(\partial
M,E_{/\partial
M})$ is self adjoint. The operator $P_{APS}$ defining the boundary
condition is the orthogonal projection onto the closed subspace of $%
L^2(\partial M,E_{/\partial M})$ spanned by the eigenfunctions of $A$
associated to non negative eigenvalues.

The projector $P_{APS}$ is a zero order pseudo differential operator
and its
principal symbol coincides with the one of the corresponding
Calder\'on
projector \cite{Booss-W}.

The problem (\ref{BoundaryProblem}) with $B=P_{APS}$ has a solution
$\varphi \in
H^1(M,E)$ for any $(\chi ,f)$ with $\chi $ in a finite codimensional
subspace of $L^2(M,E)$ and $f$ in the intersection of
$H^{1/2}(\partial
M,E_{/\partial M})$ with the image of $P_{APS}$. The solution is
unique up
to a finite dimensional kernel. Note that, since the codimension of
$%
P_{APS}\ [L^2(\partial M,E_{/\partial M})]$ is not finite, the
operator
\begin{equation}
(D,P_{APS}T):H^1(M,E)\rightarrow L^2(M,E)\times H^{1/2}(\partial
M,E_{/\partial M})
\end{equation}
is not Fredholm.

Definition 1 does not encompass Atiyah, Patodi and Singer (A.P.S.)
conditions since $P_{APS}$ takes values in $L^2(\partial
M,E_{/\partial M})$
instead of $L^2(\partial M,G),$ with $G$ a $r$ dimensional vector
bundle ($%
r=rank\ q(x;\xi ))$, as required in that definition. However, it is
possible
to define elliptic boundary problems according to Definition 1 by
using
conditions $\grave a\ la$ APS. For instance, the following
self-adjoint
boundary problem for the two-dimensional full Euclidean Dirac
operator is
elliptic:
\begin{equation}
\begin{array}{c}
\left(
\begin{array}{cc}
0 & D^{\dagger } \\
D & 0
\end{array}
\right) \left(
\QDATOP{f}{g}\right) =0\ \  \rm{ in }\ M, \\  \\
(P_{APS},\sigma (I-P_{APS})\ \sigma ^{*})\left( \QDATOP{f}{g}\right)
=0\ \ \rm{
on }\ \partial M.
\end{array}
\end{equation}

In fact, as mentioned above, the principal symbol of $P_{APS}$ is
equal to
the principal symbol of the Calder\'on projector associated to $D$.
So, from
(\ref{qb}) we get $\sigma _0(P_{APS})=H(\xi )$. By taking adjoints we
obtain
$\sigma _0(\sigma (I-P_{APS})\ \sigma ^{*})=H(-\xi )$. Then, the
principal
symbol of $B=(P_{APS},\sigma (I-P_{APS})\ \sigma ^{*})$ is \linebreak
$b(x;\xi )=(H(\xi ),H(-\xi ))$ and satisfies
\begin{equation}
rank(b(x;\xi )\ q(x;\xi ))=rank(q(x;\xi ))\qquad \forall \xi \neq 0.
\end{equation}

\section{Complex powers and regularized determinants for elliptic
boundary
problems.}

In this section we describe Seeley's construction of the complex
powers
of the operator $D$ under local elliptic boundary condition $B.$

We will denote by $\sigma ^{\prime }(D)$ the partial symbol of $D$,
i.e. the
symbol $\sigma (D)$ evaluated at $t=0$ and $\tau =-i\partial _t$ .

\bigskip
{\bf Definition} {\bf 2}:

The elliptic boundary problem (\ref{BoundaryProblem}) admits a cone
of
Agmon's directions if there is a cone $\Lambda $ in the $\lambda
$-complex
plane such that

1) $\forall \bar x\in M$, $\forall \bar \xi \neq 0,$ $\Lambda $
contains no
eigenvalues of the matrix $\sigma _1(D)(\bar x,\bar \xi )$ .

2) $\forall \xi :\vert \xi \vert  \geq 1,$ $rank(b(x;\xi )\ q(\lambda
)(x;\xi
))=rank(\ q(\lambda )(x;\xi ))$, $\forall \lambda \in \Lambda $,

\noindent where $q(\lambda )$ denotes the principal symbol of the
Calder\'on
projector $Q(\lambda )$ associated to $D-\lambda I$, with $\lambda $
included in $\sigma _1(D-\lambda I)$ (i.e. considering $\lambda $ of
degree
one in the expansion of $\sigma ($$D-\lambda I)$ in homogeneous
functions )
\cite{Gil} \cite{Seeley1}.

Condition 2 is equivalent to the following:

$\forall \lambda \in \Lambda ,$ $\forall x\in \partial M,$ $\forall
g\in
{\bf C}^{r}$, the initial value problem%
$$
\begin{array}{c}
\sigma _1^{\prime }(D)(x;\xi )\ u(t)=\lambda \ u(t) \\
\\
b(x;\xi )\ u(t)\vert _{t=0}=g
\end{array}
$$
has, for each $\xi \neq 0,$ a unique solution satisfying $\lim
\limits_{t\rightarrow \infty }\ u(t)=0$. This is the form under which
this condition is
stated in \cite{Seeley1}.

\bigskip\

An expression for $q(\lambda )(x;\xi )$ is obtained from (\ref{q1}):

\begin{equation}
\label{qlam}q(\lambda )(x;\xi )=\frac 1{2\pi i}\int_\Gamma \bigl(
a_1^{-1}(x,0;0,1;0)\ a_1(x,0;\xi ,0;\lambda )-z\bigr) ^{-1}dz,
\end{equation}
where $\ a_1(x,t;\xi ,\tau ;\lambda )=\sigma _1(D-\lambda I)$ , with
$%
\lambda $ considered of degree one as stated above.

As an example, we now compute $q(\lambda )(x;\xi )$ for the
two-dimensional
Dirac operator in (\ref{fd2}) on a disk. In polar coordinates the
principal
symbol of $D(A)-\lambda I$ is
\begin{equation}
\ a_1(\theta ,t;\xi ,\tau ;\lambda )=\tau \ \gamma _r-\xi \ \gamma
_\theta
-\lambda \ Id_{2\times 2},
\end{equation}
where\ \/
\begin{equation}
\gamma _r=\left(
\begin{array}{cc}
0 & e^{-i\theta } \\
e^{i\theta } & 0
\end{array}
\right) ,\qquad \gamma _\theta =\left(
\begin{array}{cc}
0 & -ie^{-i\theta } \\
ie^{i\theta } & 0
\end{array}
\right) .
\end{equation}
So
\begin{equation}
\begin{array}{c}
q(\lambda )(x;\xi )={\displaystyle{1 \over {2\pi i}}}
\displaystyle \int_\Gamma \Bigl( \gamma _r(-\xi \
\gamma _\theta -\lambda \ Id_{2\times 2})-z\Bigr) ^{-1}dz \\
\\
=\dfrac 1{2\sqrt{\xi ^2-\lambda ^2}}\left(
\begin{array}{cc}
\xi +\sqrt{\xi ^2-\lambda ^2} & -i\lambda e^{-i\theta } \\  &  \\
-i\lambda e^{i\theta } & -\xi +\sqrt{\xi ^2-\lambda ^2}
\end{array}
\right).
\end{array}
\end{equation}
Note that, for $\lambda =0$, (\ref{q2}) is recovered.

\bigskip\

Henceforth, we assume the existence of an Agmon's cone $\Lambda $.
Moreover,
we will consider only boundary conditions $B$ giving rise to a
discrete
spectrum $sp(D_B).$ Note that, this is always the case for elliptic
boundary
problems unless $sp(D_B)$ is the whole complex plane (see, for
instance,\cite
{Hor})$.$ Now, for $\vert \lambda \vert  $ large enough, $sp(D_B)\cap
\Lambda $
is empty, since there is no $\lambda $ in $sp(\sigma _1(D_B))\cap
\Lambda $.
Then, $sp(D_B)\cap \Lambda $ is a finite set.

For $\lambda \in \Lambda $ not in $sp(D_B)$, an asymptotic expansion
of the
symbol of $R(\lambda )=(D_B-\lambda I)^{-1}$ can be explicitly given
\cite
{Seeley1}:

\begin{equation}
\label{AE}\sigma (R(\lambda ))\sim \sum_{j=o}^\infty
c_{-1-j}-\sum_{j=o}^\infty d_{-1-j}
\end{equation}
where the {\it Seeley coefficients } I $c_{-1-j}$ and $d_{-1-j}$
satisfy

\begin{equation}
\label{7}\sum_{j=o}^\infty a_{1-j}\ \ \circ \ \sum_{j=0}^\infty
c_{-1-j}=I
\end{equation}
with $a_{1-j}$ homogeneous of degree $1-j$ in $(\bar \xi ,\lambda )$
defined
by
\begin{equation}
\label{8}\sigma (D-\lambda I)=\sum_{j=o}^\infty a_{1-j},
\end{equation}
$\circ $ denoting the usual composition of homogeneous symbols, and
\begin{equation}
\label{9}\left\{
\begin{array}{c}
\sigma ^{\prime }(D-\lambda )\ \circ \ \sum\limits_{j=o}^\infty
d_{-1-j}=0
\\
\\
\sigma ^{\prime }(B)\ \circ \ \sum\limits_{j=o}^\infty
d_{-1-j}=\sigma (B)\
\circ \ \ \sum\limits_{j=0}^\infty c_{-1-j}\quad \rm{at}\ \ t=0 \\
\\
\lim \limits_{t\rightarrow \infty }\ d_{-1-j}=0
\end{array}
\right.
\end{equation}
where the terms of $\sigma ^{\prime }(D-\lambda )\ $ are grouped
according
to their degree of homogeneity in $(\frac 1t,\xi ,\partial _t,\lambda
).$

Note that condition 2) implies the existence and unicity of the
solution of
(\ref{9}).

Written in more detail, the first line in (\ref{9}) becomes
\cite{Seeley1}
\begin{equation}
\label{11}a^{(1)}d_{-1-j}+\sum_{\QATOP{l<j}{k-\vert \alpha \vert
-1-l=-j}}
\frac{i^\alpha }{\alpha !}\frac{\partial ^\alpha }{\partial \xi
^\alpha }
a^{(k)}\frac{\partial ^\alpha }{\partial x^\alpha }d_{-1-l}=0,
\end{equation}
with
\begin{equation}
\label{11.1}a^{(j)}(x,t;\xi ,\partial _t;\lambda
)=\sum\limits_{l-k=j}\frac{
t^k}{k!}\frac{\partial ^k}{\partial t^k}a_l(x,t;\xi ,\partial
_t;\lambda
)\vert _{t=0} ,
\end{equation}
while the second one becomes
\begin{equation}
\label{12}
\begin{array}{c}
b_0d_{-1-j}+ \displaystyle \sum\limits_{_{
\QATOP{l<j}{k-\vert \alpha \vert -1-l=-j}}}\dfrac{i^\alpha }{\alpha
!}\dfrac{%
\partial ^\alpha }{\partial \xi ^\alpha }b_{-k}\dfrac{\partial
^\alpha }{%
\partial x^\alpha }d_{-1-l} \\  \\  \displaystyle
=\sum\limits_{\QATOP{l<j}{k-\vert \beta \vert -1-l=-j}}\dfrac{i^\beta
}{\beta !}
\dfrac{\partial ^\beta }{\partial \bar \xi ^\beta
}b_{-k}\dfrac{\partial
^\beta }{\partial \bar x^\beta }c_{-1-l}\vert _{t=0}
\end{array}
\end{equation}
It is worth noticing that, although
\begin{equation}
\label{Aet}\sigma (R(\lambda ))=\ \sum_{j=0}^\infty c_{-1-j}\ ,
\end{equation}
is an asymptotic expansion of $\sigma (R(\lambda )),$ the fundamental
solution of ($D_B-\lambda )$ obtained by Fourier transforming
(\ref{Aet})
does not in general satisfy the required boundary conditions. The
coefficients $d_{-1-j}$ must be added to the expansion in order to
correct
this deficiency.

The coefficients $c_{-1-j}$ $(x,t;\xi ,\tau ;\lambda )$ and $%
d_{-1-j}(x,t;\xi ,\tau ;\lambda )$ are meromorphic functions of
$\lambda $
with poles at those points where $\det [\sigma _1(D-\lambda )(x,t;\xi
,\tau
)]$ vanishes. The $c_{-1-j}$'s are homogeneous of degree $-1-j$ in (
$\xi
,\tau ,\lambda )$ ; the \linebreak $d_{-1-j}$ 's are also homogeneous
of degree $-1-j,$
but in ( $\frac 1t,\xi ,\tau ,\lambda )\cite{Seeley1}.$

Moreover, it can be proved from (\ref{AE}) that, for $\lambda \in
\Lambda ,$%
\begin{equation}
\label{SG}\parallel R(\lambda )\parallel _{L^2}\leq C\vert \lambda
\vert  ^{-1}
\end{equation}
with C a constant \cite{Seeley1,Gil}.

Estimate (\ref{SG}) allows for expressing the complex powers of $D_B$
as
\begin{equation}
\label{CP}(D_B)^z=\frac i{2\pi }\int_\Gamma \lambda ^z\ R(\lambda )\
d\lambda
\end{equation}
for $Re\ z<0$ , where $\Gamma $ is a closed path lying in $\Lambda $,
enclosing the spectrum \linebreak of $(D_B)$ \cite{Seeley2}$.$ Note
that such a curve $%
\Gamma $ always exists for $sp(D_B)\cap \Lambda $ finite.

For $Re\ z\geq 0$ , one defines%
\begin{equation}
(D_B)^z=(D)^l\circ (D_B)^{z-l}\ ,
\end{equation}
for $l$ a positive integer such that $Re\ (z-l)<0$ .

If $Re(z)<-\nu $, the power $(D_B)^z$ is an integral operator with
continuous kernel $%
J_z(x,t;y,s)$ and, consequently, it is trace class (for an
operator of order $\omega $, this is true if $Re(z)<-\frac \nu \omega
$). As
a function of $z$, $Tr(D_B)^z$ can be meromorphicaly extended to the
whole
complex plane {\bf C}, with only single poles at $z=j-\nu ,\
j=0,1,2,...$
and vanishing residues when $z=0,1,2,...$ (for an operator of order
$\omega $
, there are only single poles at $z=\frac{j-\nu }\omega ,\
j=0,1,2,...$,
with vanishing residues at $z=0,1,2,...$) \cite{Seeley2}. Throughout
this paper,  analytic functions and their meromorphic
extensions will be given the same name.

 The function $Tr(D_B)^z$ is
usually called $\zeta _{(D_B)}(-z)$ because of its similarity with
the
classical Riemann $\zeta $-function: if $\{\lambda _j\}$ are the
eigenvalues
of $D_B$, $\{\lambda _j^z\}$ are the eigenvalues of $(D_B)^z$; so $%
Tr(D_B)^z=\sum \lambda _j^z$ when $(D_B)^z$ is a trace class
operator.

A regularized determinant of $D_B$ can then be defined as
\begin{equation}
\label{DR}Det\ (D_B)=\exp [-\frac d{dz}\ Tr\ (D_B)^z]\vert _{z=0} .
\end{equation}

Now, let $D(\alpha )$ be a family of elliptic differential operators
on $M$
sharing their principal symbol
and  analytically
depending on $\alpha $. Let $B$ give rise to an elliptic boundary
condition for all of them, in such a way that $D(\alpha )_B$ is
invertible
and the boundary problems they define have a common Agmon's cone.
Then, the
variation of $Det~D(\alpha )_B$ with respect to $\alpha $ is given by
(see, for example, \cite{APS,Forman})
\begin{equation}
\label{DD}\frac d{d\alpha }\ln \ Det~D(\alpha )_B=\frac d{dz}\left[ \
z\
Tr\{\frac d{d\alpha }\left( D(\alpha )_B\right) \ (D(\alpha
)_B)^{z-1}\}\right] _{z=0} .
\end{equation}
Note that, under the assumptions made, $\frac d{d\alpha }\left(
D(\alpha )_B\right) $ is
a multiplicative operator.

Although $J_z(x,t;x,t;\alpha),$ the kernel of $(D(\alpha)_B)^z$
evaluated at the diagonal,
can be extended to the whole $z$-complex plane as a meromorphic
function, the r.h.s. in (\ref{DD}) cannot be simply written as
the integral over $M$ of the finite part of
\begin{equation}
tr\{\frac d{d\alpha }\left( D(\alpha )_B\right) \
J_{z-1}(x,t;x,t;\alpha)\}
\end{equation}
at $z=0$ (where $tr$ means matrix trace). In fact,
$J_{z-1}(x,t;x,t;\alpha)$ is in general non integrable in the
variable $t$ near $%
\partial M$ for $z\approx 0$.

Nevertheless, an integral expression for $\frac d{d\alpha }\ln \
Det~D(\alpha )_B$ will be constructed in Section 4, from the integral
expression for $Tr(D(\alpha)_B)^{z-1}$ holding in a neighborhood of
$z=0$ and obtained in the following way \cite{Seeley2}:

if $T>0$ is small enough, the function $j_z(x;\alpha)$ defined as
 \begin{equation}
\label{T}j_z(x;\alpha)\ =\int_0^TJ_z(x,t;x,t;\alpha)\ dt
\end{equation}
for $Re\ z<1-\nu $, admits a meromorphic extension to {\bf C} as a
function
of $z$. So, if $V$ is a neighborhood of $ \partial M$ defined by
$t<\epsilon $, with $\epsilon $ small enough, $Tr(D(\alpha)_B)^{z-1}$
can be written as the finite part of
\begin{equation}
\label{EI}
\int_{M/ V}tr\ J_{z-1}(x,t;x,t;\alpha)\ dxdt+\int_{\partial  M}tr\
j_{z-1}(x;\alpha)\ dx\ ,
\end{equation}
where a suitable partition of the unity is understood.


\section{Green functions and determinants}

In this section, we will give an expression for $\frac d{d\alpha }\ln
\
Det[D(\alpha )_B]$ in terms of $G_B(x,t;y,s;\alpha )$%
, the Green function of $D(\alpha )_B$ (i.e.,
the kernel of the operator $[D(\alpha )_B]^{-1}).$

With the notation of the previous Section, (\ref{DD}) can be
rewritten as:
\begin{equation}
\label{DDD}\frac d{d\alpha }\ln \ Det~D(\alpha )_B=\begin{array}{c}
\\
F.P. \\
^{_{z=0}}
\end{array}
 \int_Mtr\left[
\frac d{d\alpha }\left( D(\alpha )_B\right) \ J_{-z-1}(x,t;
x,t;\alpha
)\right] \ d\bar x\ ,
\end{equation}
where the r.h.s. must be understood as the finite part of the
meromorphic extension of the integral at $z=0$.

The finite part of $J_{-z-1}(x,t;x,t;\alpha )$ at $z=0$ does not
coincide with the regular part of $G_B(x,t;y,s;\alpha )$ at the
diagonal,
since the former is defined through an analytic extension.

However, we will show that there exists a relation between them,
involving a
finite number of Seeley's coefficients. In fact, for boundaryless
manifolds
this problem has been studied in \cite{Jour.1}, by comparing the
iterated \linebreak
limits $F.P.\lim \limits_{z\rightarrow -1}\{\lim \limits_{\bar
y\rightarrow
\bar x}J_z(x,t;y,s;\alpha)\}$ and $R.P.\lim \limits_{\bar
y\rightarrow \bar
x}\{\lim \limits_{z\rightarrow -1}J_z(x,t;y,s;\alpha)\}=$\linebreak
$R.P.\lim
\limits_{\bar y\rightarrow \bar x}G_{B}(x,t;y,s;\alpha).$

In the case of manifolds with boundary, the situation is more
involved owing
to the fact that the finite part of the extension of
$J_z(x,t;x,t;\alpha)$ \linebreak at  $z=-1$ is not integrable near
$\partial M$ $.$ (A first approach to this
problem appears in \cite{Jour.2}). Nevertheless, as mentioned in
Section 3, a meromorphic extension of $%
\int_0^T J_z(x,t;x,t;\alpha)dt,$ with $T$ small enough
can be performed and its finite part at $z=-1$ turns to be integrable
in the
tangential \linebreak variables. A similar result holds, {\it a
fortiori}, for $\int_0^T t^n J_z(x,t;x,t;\alpha)dt,$ \linebreak  with
$n=1,2,3...$
 Then, near the boundary, the Taylor expansion of the  \linebreak
function $A_\alpha =\frac d{d\alpha }D(\alpha )_B$ will naturally
appear,
and the limits \linebreak  to be compared are $F.P.\lim
\limits_{z\rightarrow -1}\{\lim
\limits_{\bar y\rightarrow \bar x}\int_0^Tt^nJ_z(x,t;y,s;\alpha)dt\}
$ and \linebreak$
R.P.\lim \limits_{\bar y\rightarrow \bar x}\{\lim
\limits_{z\rightarrow
-1}\int_0^Tt^nJ_z(x,t;y,s;\alpha)dt\}=R.P.\lim \limits_{\bar
y\rightarrow
\bar x}\int_0^Tt^nG_{B}(x,t;y,s;\alpha)dt.$

The starting point for this comparison will be to carry out
asymptotic
expansions and to analyze the terms for which the iterated limits do
not
coincide (or do not even exist).

An expansion of $G_B(x,t,y,s)$ in $M\backslash \partial M$ in
homogeneous
and logarithmic functions of $(\bar x-\bar y)$ can be obtained from (
\ref
{AE}) for $\lambda =0$:
\begin{equation}
\label{?}
\begin{array}{c}
G_B(x,t,y,s)= \sum_{j=1-\nu }^0h_j(x,t,x-y,t-s)+M(x,t)\log \vert
(x,t)-(y,s)\vert  \\
\\
+R(x,t,y,s),
\end{array}
\end{equation}
with $h_j$ the Fourier transform{\cal \ }${\cal F}^{-1}(c_{-\nu -j})$
of $%
c_{-\nu -j}$ for $j>0$ and \linebreak $h_0={\cal \ }{\cal
F}^{-1}(c_{-\nu })-\
M(x,t)\log \vert (x,t)-(y,s)\vert .$ The function $M(x,t)$ will be
explicitly
computed below (see (\ref{MM})). Our convention for the Fourier
transform is
\begin{equation}
\label{Fourier}
\begin{array}{c}
\displaystyle {\cal F}(f)(\bar \xi )=\hat f(\bar \xi )=\int f(\bar
x)\ e^{-i\bar x.\bar
\xi }\ d\bar x, \\  \\
\displaystyle {\cal \ }{\cal F}^{-1}(\hat f)(\bar x)=f(\bar x)=\dfrac
1{(2\pi )^\nu }\int
\hat f(\bar \xi )\ e^{i\bar x.\bar \xi }\ d\bar \xi .
\end{array}
\end{equation}

For $t>0$, $R(x,t,y,s)$ is continuous even at the diagonal
($y,s)=(x,t)$.
Nevertheless, $R(x,t,y,s)\vert _{(y,s)=(x,t)}$is not integrable
because of
its singularities at $t=0$. On the other hand, the functions
$t^nR(x,t,y,t)$
are integrable with respect to the variable $t$ for $y\neq x$ and $%
n=0,1,2,....$An expansion of $\int_0^\infty t^nR(x,t,y,t)dt$ in
homogeneous
and logarithmic functions of $(x-y)$ can also be obtained from
($\ref{?})$:
\begin{equation}
\label{Asn}
\int_0^\infty t^nR(x,t,y,t) dt=\sum_{j=n+2-\nu
}^0g_{j,n}(x,x-y)+M_n(x)\log (\vert x-y\vert )+R_n(x,y)
\end{equation}
where $R_n(x,y)$ is continuous even at $y=x$, and $g_{j,n}$ is the
Fourier
transform of the (homogeneous extension of) $\int_0^\infty t^n\tilde
d_{-1-j}(x,t,\xi ,t,0)dt$, with
\begin{equation}
\label{`d}\tilde d_{-1-j}(x,t,\xi ,s,\lambda )=-\int_{\Gamma
^{-}}e^{-is\tau
}\ d_{-1-j}(x,t,\xi ,\tau ,\lambda )\ d\tau
\end{equation}
for $\Gamma ^{-}$ a closed path enclosing the poles of
$d_{-1-j}(x,t,\xi
,\tau ,\lambda )$ lying in \linebreak $\{Im\ \tau >0\}$.

Since $\tilde d_{-1-j}$ is homogeneous of degree $-j$ in (1/t, $\xi
,1/s,\lambda )$, $g_{j,n}$ turns out to be homogeneous of degree $j$
in $x-y$%
{}.

The following technical lemma will be used for the proof of Theorem
1$:$

\bigskip\

{\bf Lemma 1:} {\it Let $a(\xi )$ a function defined on ${\bf R}^\nu
$,
homogeneous of \linebreak degree -$\nu $ for $\vert \xi \vert  \geq
1$ and $a(\xi )=0$ for
$\vert \xi \vert  <1.$ Then its Fourier transform can be written as
\begin{equation}
\label{4.1}{\cal \ \ }{\cal F}^{-1}(a(\xi ))(z)=h(z)+\ M\
\frac{\Omega _\nu
}{(2\pi )^\nu }\ (\log \vert z\vert ^{-1}+{\cal K}_\nu )+R(z),
\end{equation}
where

a) $h(z)$ is a homogeneous function of degree $0$, such that
\linebreak $\int_{\vert
z\vert =1}h(z)\ d\sigma _z=0.$ It is given by
\begin{equation}
\label{4.3}h(z)={\cal \ \ }{\cal F}^{-1}(P.V.[a(\xi /\vert \xi \vert
)-M\
]\vert \xi \vert  ^{-\nu })(z).
\end{equation}
b)
\begin{equation}
\label{4.2}M=\frac 1{\Omega _\nu }\int_{\vert \xi \vert  =1}a(\xi )\
d\sigma
_\xi ,
\end{equation}
where $\Omega _\nu =Area(S^{\nu -1}),$ and ${\cal K}_\nu =\ln 2-\frac
12\gamma +\frac 12\frac{\Gamma ^{\prime }(\nu /2)}{\Gamma (\nu /2)}$
with $%
\gamma $ the Euler's constant$.$

c) $R(z)$ is a function regular at $z=0$ with $R(0)=0.$ }
\bigskip

{\bf Proof:} we can decompose the function $a(\xi )$ as
\begin{equation}
\label{4.4}a(\xi )=\tilde a(\xi )+M\ \vert \xi \vert  ^{-\nu }\chi
(\xi ),
\end{equation}
where $\tilde a(\xi )$ is homogeneous of degree $-\nu $ and
$\int_{\vert \xi
\vert =1}\tilde a(\xi )\ d\sigma _\xi =0,$ and $\chi (\xi )$ is the
characteristic function of $\{\vert \xi \vert  \geq 1\}$.

Hence
\begin{equation}
\label{4.5}{\cal \ \ }{\cal F}^{-1}(\tilde a(\xi ))(z)=h(z)-r_1(z),
\end{equation}
where $h(z)$ is the Fourier transform of the distribution
\begin{equation}
\label{4.6}S=P.V.[a(\xi /\vert \xi \vert  )-M\ ]\vert \xi \vert
^{-\nu },
\end{equation}
and $r_1(z)$ is the Fourier transform of the compactly supported
distribution
\begin{equation}
\label{4.7}P.V.[\tilde a(\xi /\vert \xi \vert  )\vert \xi \vert
^{-\nu }(1-\chi
(\xi ))].
\end{equation}
Then
\begin{equation}
\label{4.8}r_1(z)=\lim _{\epsilon \rightarrow 0}\int_{\epsilon \leq
\vert \xi
\vert \leq 1}\tilde a(\xi /\vert \xi \vert  )\vert \xi \vert  ^{-\nu
}\ e^{i\xi
.z}\ {\frac{{d\xi }}{{(2\pi )^{\nu -1}}}}
\end{equation}
is a function regular at $z=0,$ with $r_1(0)=0$ . As the distribution
$S$
coincides in ${\bf R}^\nu -\{0\}$ with a smooth function of degree
$-\nu $ ,
$h={\cal F}^{-1}(S)$ is also a smooth homogeneous function of
degree $0$ and it can be seen that its mean vanishes in $S^{\nu -1}$.

On the other hand, a direct calculation gives
\begin{equation}
\label{4.9}
\begin{array}{c}
\displaystyle {\cal \ \ }{\cal F}^{-1}(\vert \xi \vert  ^{-\nu }\chi
(\xi ))(z) \\  \\
\displaystyle =
{\dfrac{{1}}{{(2\pi )^\nu }}}\int_1^\infty \dfrac{d\rho }\rho \Omega
_{\nu
-1}\int_0^\pi d\theta \ \sin {}^{\nu -2}(\theta )\ e^{i\rho \vert
z\vert \cos
\theta } \\  \\
\displaystyle =(2\pi )^{-\nu /2}\int_{\vert z\vert }^\infty d\rho \
\rho ^{-\nu /2}J_{\frac
\nu 2-1}(\rho ) \\
\\
\displaystyle =\ \dfrac{\Omega _\nu }{{(2\pi )^\nu }}\ \{\log \vert
z\vert ^{-1}+{\cal K}
_\nu \}{\sf ,}
\end{array}
\end{equation}
where
\begin{equation}
\label{4.10}
\begin{array}{c}
\displaystyle {\cal K}_\nu =\int_0^1d\rho \ \rho ^{-\nu /2}\left[
J_{\frac \nu 2-1}(\rho )-
\dfrac{\rho ^{\frac \nu 2-1}}{2^{\frac \nu 2-1}\Gamma (\nu
/2)}\right]
+\int_1^\infty d\rho \ \rho ^{-\nu /2}J_{\frac \nu 2-1}(\rho ) \\  \\
\displaystyle =\ln 2-\frac 12\gamma +\frac 12\dfrac{\Gamma ^{\prime
}(\nu /2)}{\Gamma (\nu
/2)}
\end{array}
\end{equation}
with $\gamma $ the Euler's constant ($\gamma =0.5772...)$.\qed

\bigskip\

Now, we introduce the main result of this section.

{\bf Theorem 1:}
{\it
Let $M$ be a compact connected sub manifold of ${\bf R}^\nu $ of
codimension zero
with smooth boundary $\partial M$ and $E$ a $k$-dimensional complex
vector
bundle over $M$.

Let $(D_\alpha )_{B}$ be a family of elliptic differential operators
of first order, acting on the sections of $E$, with a fixed local
boundary condition $B$ on $ \partial M$, and denote by
$J_z(x,t;x,t;\alpha)$ the meromorphic extension of the evaluation at
the diagonal of  the kernel of  $((D_\alpha )_B)^z$.

Let us assume that, for each $\alpha $,  $(D_\alpha )_B$ is
invertible, the family is differentiable with respect to $%
\alpha ,$ and $\dfrac \partial {\partial \alpha }(D_\alpha
)_{B}f=A_\alpha f$ , with $A_\alpha $ a differentiable function.

If $V$ is a neighborhood  of  $\partial M$ defined by $t < \epsilon$
and $T>0$ small enough, then:

a) }
 \nopagebreak

\begin{equation}
\label{Cuculiu}
\begin{array}{c}
\displaystyle \dfrac \partial {\partial \alpha }\ln \ Det(D_\alpha
)_B
 \\  \\ =
\begin{array}{c}
\\
F.P. \\
^{_{z=-1}}
\end{array}
\displaystyle \left[ \int_{ \partial M}\int_0^Ttr\left\{ A_\alpha
(x,t)
\  J_z(x,t;x,t;\alpha)\ \right\} dtdx\right]   \\ \\
+\begin{array}{c}
\\
F.P. \\
^{_{z=-1}}
\end{array}
\displaystyle \left[ \int_{ M/V}tr\left\{ A_\alpha (\bar x)\
J_{z}(\bar x;\bar x;\alpha)\ \right\}d\bar x\right],
\end{array}
\end{equation}
{\it where a suitable partition of the unity is understood. (This
expression must be considered as the evaluation at $z=-1$ of the
analytic extension).

b) For every $\alpha$, the integral $\int_0^TA_\alpha (x,t)\ J_z
(x,t;x,t;\alpha)dt\ $ is a meromorphic
function of $z$, for each $x\in \partial M$, with a simple pole at
$z=-1$.
Its finite part (dropping, from now on, the index $\alpha$ for the
sake of simplicity) is given by}
\begin{equation}
\label{rerptmdmadre}
\begin{array}{c}
\begin{array}{c}
\\
\displaystyle F.P. \\
^{_{z=-1}}
\end{array}
\displaystyle \int_0^TA(x,t)\ J_z(x,t;x,t)dt\\
\displaystyle = -\int_0^TA(x,t)\int_{\vert (\xi ,\tau )\vert =1}\frac
i{2\pi }\int_\Gamma
\dfrac{\ln \lambda }\lambda \ c_{-\nu }(x,t;\xi ,\tau ;\lambda )\
d\lambda \
\dfrac{d\sigma _{\xi ,\tau }}{{(2\pi )^\nu }}\ dt\  \\  \\
\displaystyle +\sum\limits_{l=0}^{\nu -2}
\dfrac{\partial _t^lA(x,0)}{l!}\int_{\vert \xi \vert
=1}\int_0^\infty t^l\
\frac i{2\pi }\int_\Gamma \dfrac{\ln \lambda }\lambda \ \tilde
d_{-(\nu
-1)+l}(x,t;\xi ,t;\lambda)\ d\lambda \ \ dt\ \dfrac{d\sigma _\xi
}{{{(2\pi )^{\nu
-1}}}}\  \\  \\
\displaystyle +\lim \limits_{y\rightarrow x}\left\{
\int_0^TA(x,t)\left[
G_B(x,t;y,t)-\sum\limits_{l=1-\nu }^0h_l(x,t;x-y,0)\right. \right.
\\
\\
\displaystyle \left. -M(x,t)\
\dfrac{\Omega _\nu }{{(2\pi )^\nu }}\left( \ln \vert x-y\vert
^{-1}+{\cal K}%
_\nu \right) \right] dt \\  \\
\displaystyle +\sum\limits_{j=0}^{^{\nu -2}}\sum\limits_{l=0}^{\nu
-2-j}
\dfrac{\partial _t^lA(x,0)}{l!}\ g_{j,l-(\nu -2-j)}(x,x-y)\  \\  \\
\displaystyle +\left. \sum_{l=0}^{^{\nu -2}}\dfrac{\partial
_t^lA(x,0)}{l!}\ M_{\nu
-2-l}(x)\ \dfrac{\Omega _{\nu -1}}{{{(2\pi )^{\nu -1}}}}\left( \ln
\vert
x-y\vert ^{-1}+{\cal K}_{\nu -1}\right) \right\} ,
\end{array}
\end{equation}
{\it
with
}
\begin{equation}
\label{MM}
\begin{array}{c}
\displaystyle M(x,t)=\dfrac 1{\Omega _\nu }\int_{\vert (\xi ,\tau
)\vert =1}c_{-\nu
}(x,t;\xi ,\tau ;0)\ \ d\sigma _{\xi ,\tau } \\
\\
\displaystyle M_j(x)=\dfrac 1{\Omega _{\nu -1}}\int_{\vert \xi \vert
=1}\int_0^\infty t^{\nu
-2-j}\ \ \tilde d_{-1-j}(x,t;\xi ,t;0)\ dt\ d\sigma _\xi ,
\end{array}
\end{equation}
{\it
where $\Omega _n=Area(S^{n-1})$, and where $h_{l}$ and $g_{l}$
are related to the Green function $G_B$ as in  (\ref{?}) and
(\ref{Asn})
}
\begin{equation}
\label{hg}
\begin{array}{c}
\displaystyle h_{1-\nu +j}(x,t;w,u)
 = {\cal F}_{(\xi ,\tau )}^{-1}\ \left[ c_{-1-j}(x,t;(\xi ,\tau
)/\vert (\xi ,\tau )\vert ;0)\ \vert (\xi ,\tau )\vert ^{-1-j}\right]
(w,u), \\ \\  \\

\displaystyle h_0(x,t;w,u)
= {\cal F}_{(\xi ,\tau )}^{-1}\left[ P.V.\left\{ \left( c_{-\nu
}(x,t;(\xi ,\tau )/\vert (\xi ,\tau )\vert ;0)-M(x,t)\right) \ \vert
(\xi ,\tau
)\vert ^{-\nu }\right\} \right] (w,u), \\  \\  \\

\displaystyle g_{j,l}(x,w)
 = {\cal F}_\xi ^{-1}\left[ \int_0^\infty t^n\ \tilde
d_{-1-j}(x,t;\xi /\vert \xi \vert  ,t;0)\ dt\vert \xi \vert
^{-1-j-n}\right] (w), \\  \\  \\

\displaystyle g_{j,0}(x,w)
={\cal F}_\xi ^{-1}\left[ P.V.\left[ \int_0^\infty t^{\nu -j-2}\
\tilde d_{-1-j}(x,t;\xi /\vert \xi \vert  ,t;0)\ dt-M_j(x)\right]
\vert \xi \vert
^{-(\nu -1)}\right] (w).
\end{array}
\end{equation} \\
{\it
c) The integral $\int_{M\backslash V}tr\left[ A_\alpha (\bar x)\
J_{z}(\bar x;\bar
x)\right] \ d\bar x $ in the second term in the r.h.s. of
(\ref{Cuculiu}) ,  is a meromorphic function of $z$ with a simple
pole at $z=-1$. Its finite part is given by
}
\begin{equation}
\begin{array}{c}
\begin{array}{c}
\\
F.P. \\
^{_{z=-1}}
\end{array}
\displaystyle \int_{M\backslash V}tr\left[ A_\alpha (\bar x)\
J_{z}(\bar x;\bar
x)\right] \ d\bar x \\
\\
\displaystyle =\int_{M\backslash V}A_\alpha (\bar x)\int_{\vert \bar
\xi \vert =1}\dfrac
i{2\pi }\int
\dfrac{\ln \lambda }\lambda \ c_{-\nu}(\bar x,\bar \xi ;\lambda )\
d\lambda \dfrac{d\bar \xi }{(2\pi )^\nu } \\  \\
\displaystyle +\int_{M\backslash V}\lim \limits_{\bar y\rightarrow
\bar x}\ A_\alpha (\bar
x)[G_B(\bar x,\bar y)-\sum\limits_{l=1-\nu }^0h_l(\bar x,\bar x-\bar
y) \\
\\
\displaystyle -M(\bar x)\dfrac{\Omega _\nu }{(2\pi )^\nu }(\ln \vert
\bar x-\bar y\vert
^{-1}+{\cal K}_\nu )]\ d\bar x.
\end{array}
\end{equation}

{\bf Proof:}
Statement  a) is a direct consequence of  (\ref{DD}), (\ref{T}) and
(\ref{EI}).

For proving b), we first establish a technical result obtained from
the fundamental estimate
\begin{equation}
\label{4.11}\vert t^n\partial _\xi ^\alpha \tilde d_{-1-j}(x,t,\xi
,s;\lambda
)\vert \leq Ce^{-c(t+s)(\vert \xi \vert  +\vert \lambda \vert
)}(\vert \xi \vert
+\vert \lambda \vert  )^{-j-n-\vert \alpha \vert },
\end{equation}
for $t,s>0$ , $\lambda \in \Lambda $, due to R.T. Seeley
\cite{Seeley1}:
 \eject
{\bf Lemma 2:}

{\it Let us define
\begin{equation}
\label{DDDD}D_{-1-j}(x,t;\xi ,t;z)\ \equiv \frac i{2\pi }\int_\Gamma
\lambda
^z\theta _1(\xi ,\lambda )\ \tilde d_{-1-j}(x,t;\xi ,t;\lambda )\
d\lambda ,
\end{equation}
then

i) If $r(x,t)$ is a function satisfying $\vert r(x,t)\vert \leq
Ct^{n}$
for $0<t<T$, $n\in {\bf N,}$ $T>0,$
\begin{equation}
\label{4.12}\int_0^Tr(x,t)\int_{{\bf R}^{\nu -1}}D_{-1-j}(x,t;\xi
,t;z)\
e^{i(x-y)\xi }\ {\frac{{d\xi}}{{(2\pi)^{\nu -1}}}} \ dt
\end{equation}
is an absolutely convergent integral for $Re(z)<j+n-\nu +1.$ As a
consequence, it is an analytic function of $z$ in this region, and it
is
continuous in all the variables ($x,y,z).$

ii) If $x\neq y$, (\ref{4.12}) is an absolutely convergent integral
for all
$z\in {\bf C,\ }$and so no analytic extension is needed out of the
diagonal.

iii)
\begin{equation}
\label{4.13}\int_0^\infty t^nD_{-1-j}(x,t;\xi ,t;z)\ \ dt
\end{equation}
is an homogeneous function of $\xi $ for $\vert \xi \vert  \geq 1$,
of degree $%
z-j-n$, analytic in $z$ for $Re(z)<j+n$ and then
\begin{equation}
\begin{array}{c}
\label{4.14}\displaystyle \int_{{\bf R}^{\nu -1}}\int_0^\infty
t^nD_{-1-j}(x,t;\xi ,t;z)\
\ dt\ {\frac{{d\xi }}{{(2\pi )^{\nu -1}}}}\\ \\ \displaystyle =\alpha
_j^n(x;z)
+\frac 1{z-j-n+\nu -1}\beta _j^n(x;z)
\end{array}
\end{equation}
with
\begin{equation}
\label{4.15}
\begin{array}{c}
\displaystyle \alpha _j^n(x;z)=\int_{\vert \xi \vert  \leq
1}\int_0^\infty
t^nD_{-1-j}(x,t;\xi ,t;z)\ \ dt\
{\dfrac{{d\xi }}{{(2\pi )^{\nu -1}}}} \\  \\
\displaystyle \beta _j^n(x;z)=\int_{\vert \xi \vert  =1}\int_0^\infty
t^nD_{-1-j}(x,t;\xi
,t;z)\ \ dt\ {\dfrac{{d\xi }}{{(2\pi )^{\nu -1}}}}
\end{array}
\end{equation}
analytic functions of $z$ for $Re(z)<j+n$.

iv)
\begin{equation}
\label{4.16}\int_{{\bf R}^{\nu -1}}\int_T^\infty t^nD_{-1-j}(x,t;\xi
,t;z)\
e^{i(x-y)\xi }\ dt\ {\frac{{d\xi}}{{(2\pi)^{\nu -1}}}}
\end{equation}
is an entire function of $z$, continuous in $(x,y,z).$ }

 \bigskip
{\bf Proof:}

i) $\tilde d_{-1-j}(x,t;\xi ,s;\lambda )$ is a homogeneous function
of
degree $-j$ in the variables $(\xi ,t^{-1},s^{-1},\lambda )$
\cite{Seeley1}%
.\ Then
\begin{equation}
\label{4.17}D_{-1-j}(x,t;\xi ,t;z)\ =\vert \xi \vert
^{z-j+1}D_{-1-j}(x,\vert
\xi \vert t;\xi /\vert \xi \vert  ,\vert \xi \vert  t;z)\
\end{equation}
for $\vert \xi \vert  \geq 1.$ In fact,
\begin{equation}
\label{4.18}
\begin{array}{c}
D_{-1-j}(x,t;\xi ,t;z)\ =\displaystyle \frac i{2\pi }\int_\Gamma
\lambda ^z\theta _1(\xi
,\lambda )\ \tilde d_{-1-j}(x,t;\xi ,t;\lambda )\ d\lambda \\
\\
=\displaystyle \vert \xi \vert  ^{-j}\frac i{2\pi }\int_\Gamma
\lambda ^z\ \tilde
d_{-1-j}(x,\vert \xi \vert  t;\xi /\vert \xi \vert  ,\vert \xi \vert
t;\lambda /\vert
\xi \vert )\ d\lambda
\end{array}
\end{equation}
because $\theta _1(\xi ,\lambda )=1$ for $\vert \xi \vert \geq 1$;
and
taking $\mu =\lambda /\vert \xi \vert  $, we get
\begin{equation}
\label{4.19}\vert \xi \vert  ^{z-j+1}\frac i{2\pi }\int_{\Gamma _\xi
}\mu ^z\
\tilde d_{-1-j}(x,\vert \xi \vert  t;\xi /\vert \xi \vert  ,\vert \xi
\vert  t;\mu )\
d\mu ,
\end{equation}
where $\Gamma _\xi =\{\lambda /\vert \xi \vert  :\lambda \in \Gamma
\}$.
Since $\tilde d_{-1-j}(x,\vert \xi \vert  t;\xi /\vert \xi \vert
,\vert \xi \vert
t;\mu )\ $has no poles between the paths $\Gamma _{\xi }$ and $
\Gamma $ for $\vert \xi \vert \geq 1$ \cite{Seeley1}, one can take $
\int_\Gamma d\lambda $ in (\ref{4.19}), thus obtaining (\ref{4.17}).

For proving i), it is sufficient to see that the integrand in
(\ref{4.12})
is dominated by an absolutely integrable function.

We write
\begin{equation}
\label{4.20}
\begin{array}{c}
\displaystyle \int_0^T\int_{
{\bf R}^{\nu -1}}\vert r(x,t)D_{-1-j}(x,t;\xi ,t;z)\vert \ \
{\dfrac{{d\xi }}{{%
(2\pi )^{\nu -1}}}}\ dt \\  \\
\displaystyle \leq C\int_0^T\int_{\vert \xi \vert  \leq 1}t^n\vert
D_{-1-j}(x,t;\xi ,t;z)\vert
\ \
{\dfrac{{d\xi }}{{(2\pi )^{\nu -1}}}}\ dt\  \\  \\
\displaystyle +C\int_0^T\int_{\vert \xi \vert  \geq 1}t^n\vert
D_{-1-j}(x,t;\xi ,t;z)\vert \ \ {
\ \dfrac{{d\xi }}{{(2\pi )^{\nu -1}}}}\ dt.
\end{array}
\end{equation}
For the first integral we have
\begin{equation}
\label{4.21}
\begin{array}{c}
\displaystyle  \int_0^T\int_{\vert \xi \vert  \leq 1}t^n\vert
D_{-1-j}(x,t;\xi ,t;z)\vert \ \
{\dfrac{{d\xi }}{{(2\pi )^{\nu -1}}}}\ dt \\  \\
\displaystyle \leq \frac 1{2\pi }\int_0^T\int_{\vert \xi \vert  \leq
1}t^n\ \int_\Gamma \vert
\lambda ^z\vert \theta _1(\xi ,\lambda )\ \vert \tilde
d_{-1-j}(x,t;\xi
,t;\lambda )\vert \ \vert d\lambda \vert {\dfrac{{d\xi }}{{(2\pi
)^{\nu -1}}}}\
dt
\end{array}
\end{equation}
and, using the estimate (\ref{4.11}), we can dominate it by
\begin{equation}
\label{4.22}C\int_\Gamma \vert \lambda ^z\vert (\vert \xi \vert
+\vert \lambda
\vert )^{-j-n}\vert d\lambda \vert \leq C\int_\Gamma \vert \lambda
\vert
^{Re(z)-j-n}\vert d\lambda \vert ,
\end{equation}
which is finite for $Re(z)<j+n-1.$

For the second integral in (\ref{4.20}) we get, from (\ref{4.17}),
\begin{equation}
\label{4.23}
\begin{array}{c}
\displaystyle \int_0^T\int_{\vert \xi \vert  \geq 1}t^n\vert
D_{-1-j}(x,t;\xi ,t;z)\vert \ \
{\dfrac{{d\xi }}{{(2\pi )^{\nu -1}}}}\ dt\\ \\
\displaystyle  =\dfrac 1{2\pi
}\int_0^T\int_{\vert \xi \vert  \geq 1}t^n\ \vert \xi \vert
^{Re(z)-j+1-n} \\  \\
\displaystyle \times \vert \int_\Gamma \lambda ^z\ (\vert \xi \vert
t)^n\ \tilde
d_{-1-j}(x,\vert \xi \vert  t;\xi /\vert \xi \vert  ,\vert \xi \vert
t;\lambda /\vert
\xi \vert )\ d\lambda \vert \ {\dfrac{{d\xi }}{{(2\pi )^{\nu -1}}}}\
dt\  .
\end{array}
\end{equation}
{}From Seeley's estimate (\ref{4.11}), it is dominated by
\begin{equation}
\label{4.24}C\int_0^T\int_{\vert \xi \vert  \geq 1}\ \vert \xi \vert
^{Re(z)-j+1-n}\int_\Gamma \vert \lambda ^z\vert \ e^{-c\vert \xi
\vert  t(1+\vert
\lambda \vert )}\ (1+\vert \lambda \vert  )^{-j-n}\vert d\lambda
\vert \ {\frac{{%
d\xi }}{{(2\pi )^{\nu -1}}}}\ dt\ .
\end{equation}
By performing the integral in t, one gets
\begin{equation}
\label{4.25}C\left( \int_{\vert \xi \vert  \geq 1}\ \vert \xi \vert
^{Re(z)-j-n}\
{\frac{{d\xi }}{{(2\pi )^{\nu -1}}}}\right) \left( \int_\Gamma \vert
\lambda
^z\vert \ \ (1+\vert \lambda \vert  )^{-1-j-n}\vert d\lambda \vert \
\right) ,
\end{equation}
which is finite for $Re(z)<j+n-\nu +1.$ \\

ii) We have, for $Re(z)<j+n-\nu +1$%
\begin{equation}
\label{4.26}
\begin{array}{c}
\displaystyle  \int_0^T\int_{
{\bf R}^{\nu -1}}r(x,t)\ D_{-1-j}(x,t;\xi ,t;z)\ e^{i(x-y)\xi }\
{\dfrac{{%
d\xi }}{{(2\pi )^{\nu -1}}}}\ dt \\  \\  \\
=\displaystyle \int_{{\bf R}^{\nu -1}}\int_0^Tr(x,t)\
D_{-1-j}(x,t;\xi ,t;z)\ e^{i(x-y)\xi
}\ \ dt\ {\dfrac{{d\xi }}{{(2\pi )^{\nu -1}}}},
\end{array}
\end{equation}
from part i).

Now, the estimate
\begin{equation}
\label{4.29}
\begin{array}{c}
\displaystyle \int_0^T\vert r(x,t)\vert \ \ \int_\Gamma \vert \lambda
^z\vert \ \vert \partial
_\xi ^\alpha [\theta _1(\xi ,\lambda )\ \tilde d_{-1-j}(x,t;\xi
,t;\lambda
)]\vert \ \vert d\lambda \vert dt\  \\
\\  \\
\displaystyle \leq C\int_0^T\ \int_\Gamma \vert \lambda ^z\vert \
t^n\sum\limits_{\alpha _1\leq \alpha }C_{\alpha _1}\vert \partial
_\xi ^{\alpha
_1}\theta _1(\xi ,\lambda )\vert \ \vert \ \partial _\xi ^{\alpha
-\alpha
_1}\tilde d_{-1-j}(x,t;\xi ,t;\lambda )\vert \ \vert d\lambda \vert
dt \\   \\  \\
 \displaystyle \leq C\int_\Gamma \vert \lambda ^z\vert \ \int_0^T\
e^{-ct(\vert \xi \vert
+\vert \lambda \vert  )}(\vert \xi \vert  +\vert \lambda \vert
)^{-j-n-\vert \alpha
\vert }\ dt\vert d\lambda \vert  \\
\\  \\
=\displaystyle C\int_\Gamma \vert \lambda ^z\vert \ (\vert \xi \vert
+\vert \lambda \vert
)^{-1-j-n-\vert \alpha \vert }\ dt\vert d\lambda \vert  \\
\\  \\
\displaystyle \leq C\ (\vert \xi \vert  +\epsilon )^{-\delta \
}\left( \int_\Gamma \vert
\lambda ^z\vert \ \vert \lambda \vert  ^{\delta -1-j-n-\vert \alpha
\vert }\ \vert
d\lambda \vert \right) ,
\end{array}
\end{equation} \\
where $\delta >0$ can be chosen such that the integral in $\lambda $
is
convergent for $Re(z)<0$, implies
\begin{equation}
\label{4.28}\partial _\xi ^\alpha \left( \int_0^Tr(x,t)\
D_{-1-j}(x,t;\xi
,t;z)\ \ dt\right)
\begin{array}{c}
\\
\rightarrow  \\
\scriptstyle \xi \rightarrow \infty
\end{array}
0\ \rm{\ for}\ Re(z)<0\ \rm{and}\ \vert \alpha \vert \geq 0.
\end{equation}
Then, for $x\neq y$, writing $e^{i(x-y)\xi }=\frac{(-1)^k}{\vert
x-y\vert ^{2k}}\triangle _\xi ^ke^{i(x-y)\xi }$, where $\triangle
_\xi $ is
the Laplacian in the $\xi $-variables, and integrating by parts in
$\xi ,$
one gets that (\ref{4.26}) becomes
\begin{equation}
\label{4.27}\frac{(-1)^k}{\vert x-y\vert ^{2k}}\int_{{\bf R}^{\nu
-1}}\triangle _\xi ^k\left( \int_0^Tr(x,t)\ D_{-1-j}(x,t;\xi ,t;z)\ \
dt\right) e^{i(x-y)\xi }\ {\frac{{d\xi }}{{(2\pi )^{\nu -1}}}}.
\end{equation}

So, the expression in (\ref{4.27}) defines a holomorphic function of
$z$
for $Re(z)$ as large as we want, taking $k$ sufficiently large, since

\begin{equation}
\label{4.30}
\begin{array}{c}
\displaystyle \int_{
{\bf R}^{\nu -1}}\int_0^T\vert r(x,t)\vert \ \vert \triangle _\xi ^k\
D_{-1-j}(x,t;\xi ,t;z)\vert \ \ dt\ {d\xi } \\  \\  \\
\displaystyle \leq C\int_{\vert \xi \vert  \geq 1}\ \int_\Gamma \vert
\lambda ^z\vert \ \vert
\triangle _\xi ^k\ \ \tilde d_{-1-j}(x,t;\xi ,t;\lambda )\vert \
\vert d\lambda \vert \
{d\xi } \\  \\  \\
\displaystyle +C\int_{\vert \xi \vert  \leq 1}\ \int_{\Gamma \cap
\{\vert \lambda \vert  \geq
1\}}\vert \lambda ^z\vert \ \vert \triangle _\xi ^k\ \ \tilde
d_{-1-j}(x,t;\xi
,t;\lambda )\vert \ \vert d\lambda \vert \
{d\xi } \\  \\  \\
\displaystyle +C\int_{\vert \xi \vert  \leq 1}\ \int_{\Gamma \cap
\{\vert \lambda \vert  \leq
1\}}\vert \lambda ^z\vert \ \vert \triangle _\xi ^k\ [\theta _1(\xi
,\lambda )\
\tilde d_{-1-j}(x,t;\xi ,t;\lambda )]\vert \ \vert d\lambda \vert \
{d\xi },
\end{array}
\end{equation} \\
where the estimate $\vert \triangle _\xi ^k\ \ \tilde
d_{-1-j}(x,t;\xi
,t;\lambda )\vert \leq C\ (\vert \xi \vert  +\vert \lambda \vert
)^{-j-2k}$
guarantees the convergence for large $Re(z)$ of the first two terms,
and the last
one is convergent for all $z.$ \\

iii) The integral $\int_0^Tt^n\ D_{-1-j}(x,t;\xi ,t;z)\ \ dt$ is
absolutely
convergent for $Re(z)<j+n$ , as a consequence of the estimate
(\ref{4.11}).
Its homogeneity is obvious, and then
\begin{equation}
\label{4.31}
\begin{array}{c}
\displaystyle \int_{
{\bf R}^{\nu -1}}\int_0^\infty t^n\ D_{-1-j}(x,t;\xi ,t;z)\ \ dt\
{\dfrac{{%
d\xi }}{{(2\pi )^{\nu -1}}}} \\  \\
\displaystyle =\int_{\vert \xi \vert  \leq 1}\int_0^\infty t^n\
D_{-1-j}(x,t;\xi ,t;z)\ \ dt\
{\dfrac{{d\xi }}{{(2\pi )^{\nu -1}}}} \\  \\
\displaystyle +\int_{\vert \xi \vert  =1}\left( \int_1^\infty r^{\nu
-2+z-j-n}\ dr\right)
\left( \int_0^\infty t^n\ D_{-1-j}(x,t;\xi ,t;z)\ \ dt\right) \
\dfrac{d\sigma _\xi }{(2\pi )^{\nu -1}} \\  \\
\displaystyle =\alpha _j^n(x;z)+\frac 1{z-j-n+\nu -1}\beta _j^n(x;z),
\end{array}
\end{equation}
for $Re(z)<j+n-\nu +1.$
\eject
It is easy to verify the analyticity of $\alpha _j^n(x;z)$ and $\beta
_j^n(x;z)$ for $Re(z)<j+n$ by means of the estimate of $t^n\tilde
d_{-1-j}(x,t;\xi
,t;\lambda ).$

 iv) It is sufficient to prove that the integral is absolutely
convergent for
all $z.$ In fact,
\begin{equation}
\label{4.32}
\begin{array}{c}
\displaystyle \int_{
{\bf R}^{\nu -1}}\int_T^\infty t^n\ \int_\Gamma \vert \lambda ^z\vert
\theta
_1(\xi ,\lambda )\ \vert \tilde d_{-1-j}(x,t;\xi ,t;\lambda )\vert \
\vert d  \lambda\vert dt\ {\dfrac{{d\xi }}{{(2\pi )^{\nu -1}}}}\\ \\
\leq
\displaystyle
\int_{{\bf R}^{\nu -1}}\int_\Gamma \vert \lambda ^z\vert \left(
\int_T^\infty \ e^{-ct(\vert \xi \vert  +\vert \lambda \vert  )}\
dt\right) (\vert
\xi \vert +\vert \lambda \vert  )^{-j-n}\ \theta _1(\xi ,\lambda )
\vert d  \lambda\vert dt\ {\dfrac{{d\xi }}{{(2\pi )^{\nu -1}}}}\ .
\end{array}
\end{equation}
The integral in t is $\dfrac{\ e^{-cT(\vert \xi \vert  +\vert \lambda
\vert  )}}{%
c(\vert \xi \vert  +\vert \lambda \vert  )}$ and $(\vert \xi \vert
+\vert \lambda
\vert )^{-1-j-n}\leq C$ in the support of $\theta _1(\xi ,\lambda )$.
Then, we
can bound (\ref{4.32}) by
\begin{equation}
\label{4.33}C\left( \int_{{\bf R}^{\nu -1}}\ e^{-cT\vert \xi \vert
}{\frac{{%
d\xi }}{{(2\pi )^{\nu -1}}}}\right) \left( \int_\Gamma \vert \lambda
^z\vert \
e^{-cT\vert \lambda \vert  }\vert d\lambda \vert \right) \
\end{equation}
which is finite for all values of $z$. \qed

\bigskip\

In what follows, we will  proof the assertion in (\ref{rptmd}).
As a byproduct we will also obtain the following expression for the
residue
 \begin{equation}
\label{res}
\begin{array}{c}
\begin{array}{c}
\\
\displaystyle Res \\
^{_{z=-1}}
\end{array}
\int_0^TA(x,t)\ J_z(x,t;x,t)dt \\
\\
\displaystyle =-\int_0^TA(x,t)\int_{\vert (\xi ,\tau )\vert =1}\
c_{-\nu }(x,t;\xi ,\tau
;0)\
\dfrac{d\sigma _{\xi ,\tau }}{{(2\pi )^\nu }}\ dt\  \\  \\
\displaystyle +\sum\limits_{l=0}^{\nu -2}\dfrac{\partial
_t^lA(x,0)}{l!}\int_{\vert \xi
\vert =1}\ \int_0^\infty t^l\ \tilde d_{-(\nu -1)+l}(x,t;\xi ,t;0)\
dt\
\dfrac{d\sigma _\xi }{{{(2\pi )^{\nu -1}}}}\ .\
\end{array}
\end{equation}


We can use , as an approximation to $(D_B-\lambda )^{-1}$
\cite{Seeley1},
the parametrix

\begin{equation}
\label{100}P_K(\lambda )=\sum_\varphi \psi \left[
\sum_{j=0}^KOp(\theta _2\
c_{-1-j})-\sum_{j=0}^KOp^{\prime }(\theta _{1\ }d_{-1-j})\right] \
\varphi ,
\end{equation}
where $\varphi $ is a partition of the unity, $\psi \equiv 1$ in $%
Supp(\varphi )$,  \\
\begin{equation}
\begin{array}{c}
\theta _2(\xi ,\tau ,\lambda )=\chi (\vert \xi \vert  ^2+\vert \tau
\vert ^2+\vert
\lambda \vert ^2) \\   \\
\\
\theta _1(\xi ,\lambda )=\chi (\vert \xi \vert  ^2+\vert \lambda
\vert  ^2),
\end{array}
\end{equation} \\
with \cite{Seeley1} \\
\begin{equation}
\chi (t)=\left\{
\begin{array}{cc}
0 & t\leq 1/2 \\
1 & t\geq 1
\end{array}
\right. ,
\end{equation} \\
and  \\
\begin{equation}
\begin{array}{c}
\displaystyle Op(\sigma )h(x,t)=\int \sigma (x,t;\xi ,\tau )\ \hat
h(\xi ,\tau )\
e^{i(x\xi +t\tau )}
\dfrac{d\xi }{(2\pi )^{\nu -1}}\ \dfrac{d\tau }{2\pi }, \\  \\  \\
\displaystyle Op^{\prime }(\sigma )h(x,t)=\int \int \tilde \sigma
(x,t;\xi ,s)\ \tilde
h(\xi ,s)\ e^{ix\xi }\dfrac{d\xi }{(2\pi )^{\nu -1}}\ \dfrac{ds}{2\pi
},
\end{array}
\end{equation} \\
where $\ \hat h(\xi ,\tau )$ is defined in (\ref{Fourier}) and
\begin{equation}
\tilde h(\xi ,s)=\int \ h(x,s)\ e^{-ix\xi }\ dx.\
\end{equation} \\
Thus, one can approximate the kernel $J_z$ of $D_B^z$ by means of the
kernel
$L_z^K$ of $\frac i{2\pi }\int_\Gamma \lambda ^z\ P_K(\lambda )\
d\lambda $.
We have \\  \\
\begin{equation}
\begin{array}{c}
\displaystyle L_z^K(x,t;y,s)=\sum\limits_\varphi \ \psi (x,t)\left[ \
\sum\limits_{j=0}^K\int_{\bf{R}^\nu }C_{-1-j}(x,t;\xi ,\tau ;z)\
e^{i[(x-y)\xi
+(t-s)\tau ]}\
{\dfrac{{d\xi }}{{(2\pi )^{\nu -1}}}}\ {\dfrac{{d\tau }}{{2\pi
}}}\right.
\\  \\  \\
\displaystyle \left. -\sum\limits_{j=0}^K\int_{\bf{R}^{\nu
-1}}D_{-1-j}(x,t;\xi ,\tau ;z)\
e^{i(x-y)\xi \ }{\dfrac{{d\xi }}{{(2\pi )^{\nu -1}}}}\right] \
\varphi (y,s)
\end{array}
\end{equation}
\eject
where we have called \\
\begin{equation}
C_{-1-j}(x,t;\xi ,\tau ;z)=\frac i{2\pi }\int_\Gamma \lambda ^z\
\theta
_2(\xi ,\tau ;\lambda )\ \ c_{-1-j}(x,t;\xi ,\tau ;\lambda )\
d\lambda ,
\end{equation}
and  \\
\begin{equation}
D_{-1-j}(x,t;\xi ,t;z)=\frac i{2\pi }\int_\Gamma \lambda ^z\ \theta
_1(\xi
;\lambda )\ \ \ \tilde d_{-1-j}(x,t;\xi ,t;\lambda )\ d\lambda ,
\end{equation} \\
as in (\ref{DDDD}). These  expressions are, in fact, analytic
functions of $z$
for all complex $z,$ since the singularities of $\ c_{-1-j}(\lambda
)\ $and $%
\tilde d_{-1-j}(\lambda )\ $ are in a compact set in the $\lambda $
plane,
for $(x,t;\xi ,\tau )$ in a compact set.

Since $(D_B-\lambda )^{-1}-P_K(\lambda )$ has a continuous kernel of
$O(\
\vert \lambda \vert  ^{\nu -K-1})$ for $\lambda \in \Lambda $
\cite{Seeley1},
it turns out that\\
\begin{equation}
R(x,t;y,s;z)=J_z(x,t;y,s)-L_z^K(x,t;y,s)
\end{equation}\\
is a continuous function of $x,t,y,s$ and $z$, and analytic in $z$
for $%
Re(z)<0$, if $K\geq \nu $. Analyzing the last terms in $L_z^K$, we
obtain
that it is also true for $K=\nu -1.$ From now on, we call
$L_z=L_z^{\nu -1}.$
Then \\  \\
\begin{equation}
\lim _{z\rightarrow -1}\left[ \lim _{(y,s)\rightarrow
(x,t)}(J_z-L_z)\right]
=\lim _{(y,s)\rightarrow (x,t)}\left[ \lim _{z\rightarrow
-1}(J_z-L_z)\right] .
\end{equation}\\
Since \\
\begin{equation}
J_{-1}(x,t;y,s)=G_B(x,t;y,s),\ \  \hbox {for}\ (x,t)\neq (y,s),
\end{equation}
we have  \\  \\
\begin{equation}
\label{jl}\lim _{z\rightarrow -1}(J_z(x,t;x,t)-L_z(x,t;x,t))=\lim
_{(y,s)\rightarrow (x,t)}(G_B(x,t;y,s)-L_{-1}(x,t;y,s)). \\
\end{equation}

\bigskip
As will be shown in Lemma 3 below, one can cancel some terms in the
equality
(\ref{jl})  by studying the singularities of $L_{-1}(x,t;y,s)$ at
$(x,t)=(y,s)$, and
those of $L_z(x,t;x,t)$ at $z=-1$.

\eject
\bigskip
{\bf Lemma 3} :\\ \\ {\it The following statement holds} \\

\begin{equation}
\label{lemma3}
\begin{array}{c}
\displaystyle \lim \limits_{z\rightarrow -1}\left[ J_z(x,t;x,t)+\frac
1{(z+1)}\int_{\vert
(\xi ,\tau )\vert =1}c_{-\nu }(x,t;\xi ,\tau ;0)\
\dfrac{d\sigma _{\xi ,\tau }}{(2\pi )^\nu }\right.  \\  \\  \\ \\
\displaystyle +\int_{\vert (\xi ,\tau )\vert =1}\frac i{2\pi
}\int_\Gamma
\frac{\ln \lambda }\lambda c_{-\nu }(x,t;\xi ,\tau ;\lambda )\
d\lambda \
\dfrac{d\sigma _{\xi ,\tau }}{(2\pi )^\nu } \\  \\  \\ \\
\displaystyle \left. +\sum\limits_{j=0}^{\nu -1}\int_{
{\bf R}^{\nu -1}}D_{-1-j}(x,t;\xi ,t;z)\ \ {\dfrac{{d\xi }}{{(2\pi
)^{\nu -1}%
}}}\right]  \\  \\  \\ \\
\displaystyle =\lim \limits_{y\rightarrow x}\left[
G_B(x,t;y,t)-\sum\limits_{l=1-\nu
}^0h_l(x,t;x-y,0)\right.  \\   \\  \\ \\
\displaystyle -\ M(x,t)\
\dfrac{\Omega _\nu }{(2\pi )^\nu }(\ln \vert x-y\vert ^{-1}+{\cal
K}_\nu ) \\  \\ \\  \\
\displaystyle \left. +\sum\limits_{j=0}^{\nu -1}\int_{{\bf R}^{\nu
-1}}D_{-1-j}(x,t;\xi
,t;-1)e^{i(x-y)\xi }\ \ {\dfrac{{d\xi }}{{(2\pi )^{\nu -1}}}}\right]
{}.
\end{array}
\end{equation}
\eject

{\bf Proof:}

The terms in $L_z(x,t;x,t)$ involving the coefficients $C_{-1-j}$ can
be
written as:  \\ \\
\begin{equation}
\label{siete}
\begin{array}{c}
\displaystyle \sum\limits_{j=0}^{\nu -1}\int C_{-1-j}(x,t;\xi ,\tau
;z)\ \
{\dfrac{{d\xi }}{{(2\pi )^{\nu -1}}}}\ {\dfrac{{d\tau }}{{2\pi }}} \\
 \\  \\
\displaystyle =\sum\limits_{j=0}^{\nu -2}\int_{\vert (\xi ,\tau
)\vert \leq
1}C_{-1-j}(x,t;\xi ,\tau ;z)\ \
{\dfrac{{d\xi }}{{(2\pi )^{\nu -1}}}}\ {\dfrac{{d\tau }}{{2\pi }}} \\
 \\  \\
\displaystyle +
\dfrac{-1}{z-j-\nu }\int_{\vert (\xi ,\tau )\vert   =
1}C_{-1-j}(x,t;\xi
,\tau ;z)\ d\sigma _{\xi ,\tau } \\  \\  \\
\displaystyle +\int_{\vert (\xi ,\tau )\vert \leq 1}C_{-\nu }(x,t;\xi
,\tau ;z)\ \
{\dfrac{{d\xi }}{{(2\pi )^{\nu -1}}}}\ {\dfrac{{d\tau }}{{2\pi }}} \\
 \\  \\
\displaystyle +
\dfrac{-1}{(\ z+1)}\int_{\vert (\xi ,\tau )\vert =1}C_{-\nu }(x,t;\xi
,\tau
;z)\vert _{z=-1}\ \dfrac{d\sigma _{\xi ,\tau }}{{(2\pi )^\nu }} \\
\\  \\
\displaystyle -i\int_{\vert (\xi ,\tau )\vert =1}\partial _zC_{-\nu
}(x,t;\xi ,\tau ;z)\vert
_{z=-1}\
\dfrac{d\sigma _{\xi ,\tau }}{{(2\pi )^\nu }} \\ \\ \\
+O(z+1), \\ \\ \\
\end{array}
\end{equation}
because of the homogeneity properties of the functions $C_{-1-j}$ for
$\vert
(\xi ,\tau )\vert \geq 1$, and their analyticity in the variable $z$.

Analogously, considering $L_{-1}(x,t;x,s)$ we get
\begin{equation}
\label{ocho}
\begin{array}{c}
\displaystyle \sum\limits_{j=0}^{\nu -1}\int C_{-1-j}(x,t;\xi ,\tau
;-1)\ e^{i[(x-y)\xi
+(t-s)\tau ]}\
{\dfrac{{d\xi }}{{(2\pi )^{\nu -1}}}}\ {\dfrac{{d\tau }}{{2\pi }}} \\
 \\
\displaystyle =\sum\limits_{j=0}^{\nu -2}\Biggl\{ h_{1-\nu
+j}(x,t;x-y,t-s)\Biggr.  \\  \\
\displaystyle +\int_{\vert (\xi ,\tau )\vert \leq 1}C_{-1-j}(x,t;\xi
,\tau ;-1)\ \
e^{i[(x-y)\xi +(t-s)\tau ]}\
{\dfrac{{d\xi }}{{(2\pi )^{\nu -1}}}}\ {\dfrac{{d\tau }}{{2\pi }}} \\
 \\  \\
\displaystyle \left. -\int_{\vert (\xi ,\tau )\vert \leq
1}C_{-1-j}(x,t;(\xi ,\tau )/\vert
(\xi ,\tau )\vert ;-1)\ \ \vert (\xi ,\tau )\vert ^{-1-j}\
{\dfrac{{d\xi }}{{(2\pi )^{\nu -1}}}}\ {\dfrac{{d\tau }}{{2\pi
}}}\right\}
\\  \\  \\
\displaystyle +h_0(x,t;x-y,t-s)+\ M(x,t)
\dfrac{\Omega _\nu }{(2\pi )^\nu }\left[ \ln (\vert x-y\vert ^2+\vert
t-s\vert
^2)^{-1/2}+{\cal K}_\nu \right]  \\  \\  \\
\displaystyle +\int_{\vert (\xi ,\tau )\vert \leq 1}C_{-\nu
}(x,t;(\xi ,\tau )/\vert (\xi
,\tau )\vert ;-1)\ \ e^{i[(x-y)\xi +(t-s)\tau ]}\
{\dfrac{{d\xi }}{{(2\pi )^{\nu -1}}}}\ {\dfrac{{d\tau }}{{2\pi }}} \\
 \\  \\
+O((x-y,t-s)),
\end{array}
\end{equation}
where $h_l$,  $l>0$, are the homogeneous functions obtained by
Fourier
transforming  $C_{-1-j}(x,t;\frac{(\xi ,\tau )}{\vert (\xi
,\tau )\vert };-1)\ \vert (\xi ,\tau )\vert ^{-1-j}$.
Lemma 1 was applied for calculating the Fourier transform of
$C_{-1-j}(x,t;\xi ,\tau
;-1).$

Then, we can write (\ref{siete}) as
\begin{equation}
\begin{array}{c}
\displaystyle \dfrac{-1}{(\ z+1)}\int_{\vert (\xi ,\tau )\vert
=1}C_{-\nu }(x,t;\xi ,\tau
;-1)\ \ \dfrac{d\sigma _{\xi ,\tau }}{(2\pi )^\nu } \\  \\
\displaystyle -\int_{\vert (\xi ,\tau )\vert =1}\partial _zC_{-\nu
}(x,t;\xi ,\tau ;-1)\
\dfrac{d\sigma _{\xi ,\tau }}{(2\pi )^\nu }+R_1(x,t,z)+O(z+1)
\end{array}
\end{equation}
and expression (\ref{ocho}) as
\begin{equation}
\begin{array}{c}
\displaystyle \sum\limits_{j=0}^{\nu -1}h_{-(\nu -1)+j}(x,t;x-y,t-s)
\\
\\
\displaystyle +\ M(x,t)
\dfrac{\Omega _\nu }{(2\pi )^\nu }\left[ \ln (\vert x-y\vert ^2+\vert
t-s\vert
^2)^{-1/2}+{\cal K}_\nu \right]  \\  \\
\displaystyle +R_2(x,t,y,s)+O((x-y,t-s)) .
\end{array}
\end{equation}
Taking into account that
\begin{equation}
\lim _{z\rightarrow -1}R_1(x,t,z)=\lim _{(y,s)\rightarrow
(x,t)}R_2(x,t,y,s) ,
\end{equation}
and that
\begin{equation}
\begin{array}{c}
\displaystyle C_{-\nu }(x,t;\xi ,\tau ;-1)=c_{-\nu }(x,t;\xi ,\tau
;0)\ \
\rm{for}\ \ \vert (\xi ,\tau )\vert \geq 1 \\  \\
\displaystyle \partial _zC_{-\nu }(x,t;\xi ,\tau ;-1)=\dfrac i{2\pi
}\int_\Gamma \dfrac{%
\ln \lambda }\lambda \ c_{-\nu }(x,t;\xi ,\tau ;\lambda )\ d\lambda ,
\end{array}
\end{equation}
we obtain Lemma 3. \qed

\bigskip\

The meromorphic extension of the terms involving the coefficients
\linebreak $ C_{-1-j}(x,t;\xi ,\tau ;z)$ in $L_z(x,t;x,t)$ is a
consequence of expression
(\ref{ocho}). Although  $\sum\limits_{j=0}^{\nu -1}\int_{{\bf R}^{\nu
-1}}D_{-1-j}(x,t;\xi ,t;z)\ \ {\frac{{d\xi }}{{(2\pi )^{\nu -1}}}}$
does not
admit, in general, a meromorphic extension, such extension can be
performed
for
\begin{equation}
\int_0^Tt^n\int_{{\bf R}^{\nu -1}}D_{-1-j}(x,t;\xi ,\tau ;z)\ \
{\frac{{d\xi
}}{{(2\pi )^{\nu -1}}}}\ dt,
\end{equation}
for $n=0,1,...$ and $j=0,1,2,...$ (see \cite{Seeley2} and Lemma 2).

In order to get part b) of the theorem, we study the limits
\linebreak $\lim
\limits_{z\rightarrow -1}\int_0^TA(x,t)\ R(x,t;z)\ dt$ and $\lim
\limits_{y\rightarrow x}\int_0^TA(x,t)\ S(x,t;y,t)\ dt$, where
$R(x,t;z)$
and $S(x,t;y,t)$ denote the expressions appearing in the limits on
the
l.h.s. and r.h.s. of (\ref{lemma3}) respectively. (We have written $%
A(x,t)=A_\alpha (x,t)$ for notational simplicity).

{\bf Lemma 4: }{\it If $A(x,t)$ has $\nu -1-j$ continuous derivatives
in the
variable $t,$ $t\geq 0,$ then

i) For $\nu -1-j>0$,

\begin{equation}
\begin{array}{c}
\displaystyle \int_0^TA(x,t)\int D_{-1-j}(x,t;\xi ,t;z)\
{\dfrac{{d\xi }}{{(2\pi )^{\nu -1}}}}\ dt=\psi _j(x,z)\  \\  \\
\displaystyle -\dfrac 1{z+1}
\dfrac{\partial _t^{\nu -j-2}A(x,0)}{(\nu -j-2)!}\ \int_{\vert \xi
\vert
=1}\int_0^\infty t^{\nu -j-2}D_{-1-j}(x,t;\xi ,t;-1)\ dt\
\dfrac{d\sigma
_\xi }{{{(2\pi )^{\nu -1}}}} \\  \\
\displaystyle \ -\dfrac{\partial _t^{\nu -j-2}A(x,0)\ }{(\nu
-j-2)!}\int_{\vert \xi \vert
=1}\int_0^\infty t^{\nu -j-2}\partial _zD_{-1-j}(x,t;\xi ,t;-1)\
dt\dfrac{%
d\sigma _\xi }{{{(2\pi )^{\nu -1}}}},
\end{array}
\end{equation}
with $\psi _j(x,z)$ an analytic function of $z$ for $Re(z)<0.$

Moreover,

\begin{equation}
\begin{array}{c}
\displaystyle \int_0^TA(x,t)\int D_{-1-j}(x,t;\xi ,t;-1)\
e^{i(x-y)\xi }\
{\dfrac{{d\xi }}{{(2\pi )^{\nu -1}}}}\ dt\\  \\
\displaystyle =\varphi _j(x,y)\  +\sum\limits_{
\QATOP{n=0}{l=j+n-\nu +2}}^{\nu -j-2}\ \dfrac{\partial
_t^nA(x,0)}{n!}\
g_{j,l}(x,x-y)\  \\  \\
\displaystyle +\dfrac{\partial _t^{\nu -j-2}A(x,0)\ }{(\nu -j-2)!}\ \
M_j(x)\ \dfrac{%
\Omega _{\nu -1}}{{{(2\pi )^{\nu -1}}}}(\ln \vert x-y\vert
^{-1}+{\cal K}_{\nu
-1}),
\end{array}
\end{equation}
where $\varphi _j(x,y)$ is a continuous function.

ii) For $\nu -1-j=0$,
\begin{equation}
\int_0^TA(x,t)\int D_{-1-j}(x,t;\xi ,t;z)\
{\frac{{d\xi}}{{(2\pi)^{\nu -1}}}}
\ dt=\psi _j(x,z)
\end{equation}
is an analytic function of $z$ for $Re(z)<0$, and
\begin{equation}
\int_0^TA(x,t)\int D_{-1-j}(x,t;\xi ,t;-1)\ e^{i(x-y)\xi }\
{\dfrac{{d\xi }}{{(2\pi )^{\nu -1}}}}\ dt =\varphi _j(x,y)
\end{equation}
is a continuous function.

iii) For all $j$ }
\begin{equation}
\label{510}\lim _{z\rightarrow -1}\ \psi _j(x,z)=\lim _{y\rightarrow
x}\
\varphi _j(x,y)
\end{equation}
\bigskip

{\bf Proof:}

For $\nu -1-j\geq 0$, let
\begin{equation}
A(x,t)=\sum_{n=0}^{\nu -j-2}\frac{\partial _t^nA(x,0)\ }{n!}\ t^n\ +\
\epsilon _j(x,t)\ t^{\nu -1-j},
\end{equation}
with $\vert \epsilon _j(x,t)\vert \leq C$, for$\ t\in [0,T].$ Then,
we write \\
\begin{equation}
\label{511}
\begin{array}{c}
\displaystyle \int_0^TA(x,t)\int D_{-1-j}(x,t;\xi ,t;z)\
{\dfrac{{d\xi }}{{(2\pi )^{\nu -1}}}}\ dt \\  \\  \\
\displaystyle =\int \left( \int_0^TA(x,t)D_{-1-j}(x,t;\xi ,t;z)\ \
dt\right) \
{\dfrac{{d\xi }}{{(2\pi )^{\nu -1}}}} \\  \\  \\
\displaystyle =\sum\limits_{n=0}^{\nu -j-2}
\dfrac{\partial _t^nA(x,0)}{n!}\ \int \left( \int_0^Tt^n\
D_{-1-j}(x,t;\xi
,t;z)\ dt\right) {\dfrac{{d\xi }}{{(2\pi )^{\nu -1}}}}\\ \\  \\
\displaystyle +\int \left( \int_0^T\epsilon _j(x,t)\ t^{\nu -1-j}\
D_{-1-j}(x,t;\xi ,t;z)\
dt\right) {\dfrac{{d\xi }}{{(2\pi )^{\nu -1}}}}\ , \\ \\
\end{array}
\end{equation}
where the first equality holds for $Re(z)<j-\nu +1$ by the estimates
in
Lemma 2 i). Also, by Lemma 2 i) the last term is analytic for
$Re(z)<0$ and
it can be written as
\begin{equation}
\int \int_0^T\epsilon _j(x,t)\ t^{\nu -1-j}\ D_{-1-j}(x,t;\xi ,t;-1)\
dt\ {\
\frac{{d\xi }}{{(2\pi )^{\nu -1}}}}\ +\ O(z+1).
\end{equation}
Finally, by Lemma 2 iv), and evaluating this analytic functions at
$z=-1,$
expression (\ref{511} ) gives
\begin{equation}
\label{412}
\begin{array}{c}
\displaystyle \sum\limits_{n=0}^{\nu -j-3}
\dfrac{\partial _t^nA(x,0)}{n!}\ \left[ \alpha _j^n(x;-1)+
\dfrac{-1}{z-j-n+\nu -1}\beta _j^n(x;-1)\right.  \\  \\
\displaystyle \left. -\int \int_T^\infty t^n\ D_{-1-j}(x,t;\xi
,t;-1)\ dt\
{\dfrac{{d\xi }}{{(2\pi )^{\nu -1}}}}\right]  \\  \\
\displaystyle +
\dfrac{\partial _t^{\nu -j-2}A(x,0)}{(\nu -j-2)!}\ \left[ \alpha
_j^{\nu
-j-2}(x;-1)+\dfrac{-1}{(z+1)}\beta _j^{\nu -j-2}(x;-1)-\partial
_z\beta
_j^{\nu -j-2}(x;-1)\right.  \\  \\
\displaystyle +\left. \int \int_T^\infty t^{\nu -j-2}\ D_{-\nu
+1}(x,t;\xi ,t;-1)\ dt\
{\dfrac{{d\xi }}{{(2\pi )^{\nu -1}}}}\right]  \\  \\
\displaystyle +\int \left( \int_0^T\epsilon _j(x,t)\ t^{\nu -j-2}\
D_{-1-j}(x,t;\xi
,t;-1)\ dt\right) {\dfrac{{d\xi }}{{(2\pi )^{\nu -1}}}}+O(z+1).
\end{array}
\end{equation}
On the other hand,
\begin{equation}
\label{13}
\begin{array}{c}
\displaystyle \int_0^TA(x,t)\ \ \left( \int D_{-1-j}(x,t;\xi ,t;-1)\
e^{i(x-y)\xi }
{\dfrac{{d\xi }}{{(2\pi )^{\nu -1}}}}\right) \ dt \\  \\
\displaystyle =\int \left( \int_0^TA(x,t)\ \ \int D_{-1-j}(x,t;\xi
,t;-1)\ dt\right) \
e^{i(x-y)\xi }{\dfrac{{d\xi }}{{(2\pi )^{\nu -1}}}}
\end{array}
\end{equation}
since, for $x\neq y$ , the integral is absolutely convergent ( Lemma
2 ii)).
Then, it can be written as
\begin{equation}
\label{14}
\begin{array}{c}
\displaystyle \sum\limits_{n=0}^{\nu -j-2}
\dfrac{\partial _t^nA(x,0)\ }{n!}\int \left( \int_0^\infty t^n\
D_{-1-j}(x,t;\xi ,t;-1)\ dt\right) \ e^{i(x-y)\xi }{\dfrac{{d\xi
}}{{(2\pi
)^{\nu -1}}}} \\  \\
\displaystyle -\sum\limits_{n=0}^{\nu -j-2}
\dfrac{\partial _t^nA(x,0)\ }{n!}\int \left( \int_T^\infty t^n\
D_{-1-j}(x,t;\xi ,t;-1)\ dt\right) \ e^{i(x-y)\xi }{\dfrac{{d\xi
}}{{(2\pi
)^{\nu -1}}}} \\  \\
\displaystyle +\int \left( \int_0^T\epsilon _j(x,t)\ t^{\nu -1-j}\
D_{-1-j}(x,t;\xi
,t;-1)dt\right) \ e^{i(x-y)\xi }{\dfrac{{d\xi }}{{(2\pi )^{\nu
-1}}}}.
\end{array}
\end{equation}
Here, the integral $\int_0^\infty t^n\ D_{-1-j}(x,t;\xi ,t;-1)\ dt$
is a
homogeneous function of degree $-1-j-n$ for $\vert \xi \vert \geq 1.$
So, its
Fourier transform evaluated at $x-y$ can, for $-1-j-n\geq -(\nu -2)$,
be
written as
\begin{equation}
\label{15}
\begin{array}{c}
\displaystyle g_{j,l}(x,x-y)+\int_{\vert \xi \vert \leq 1}\left(
\int_0^\infty t^n\
D_{-1-j}(x,t;\xi ,t;-1)\ dt\ \right) e^{i(x-y)\xi }\
{\dfrac{{d\xi }}{{(2\pi )^{\nu -1}}}} \\  \\
\displaystyle -\int_{\vert \xi \vert \leq 1}\left( \int_0^\infty t^n\
D_{-1-j}(x,t;\xi /\vert
\xi \vert ,t;-1)\ dt\right) \vert \xi \vert ^{-1-j-n}e^{i(x-y)\xi }
{\dfrac{{d\xi }}{{(2\pi )^{\nu -1}}}} \\  \\
\displaystyle =g_{j,l}(x,x-y)+\alpha _j^n(x,-1)-\dfrac 1{-j-n+\nu
-2}\beta
_j^n(x,-1)+O(\vert x-y\vert ),
\end{array}
\end{equation}
where $g_{j,l}(x,w)$ is the homogeneous function of degree $l=j+n-\nu
+2$
defined in ($\ref{hg}).$

When $-1-j-n=-(\nu -1),$ by Lemma 1, we have that the Fourier
transform of $%
\int_0^\infty t^{\nu -j-2}\ D_{-1-j}(x,t;\xi ,t;-1)\ dt$ is given by

\begin{equation}
\label{16}
\begin{array}{c}
\displaystyle g_{j,0}(x,x-y)+\ M_j(x)\
\dfrac{\Omega _{\nu -1}}{{{(2\pi )^{\nu -1}}}}(\ln \vert x-y\vert
^{-1}+{\cal K%
}_{\nu -1})+O(\vert x-y\vert ) \\  \\  \\
\displaystyle +\int_{\vert \xi \vert \leq 1}\left( \int_0^\infty
t^{\nu -j-2}\
D_{-1-j}(x,t;\xi ,t;-1)\ dt\right) e^{i(x-y)\xi }
{\dfrac{{d\xi }}{{(2\pi )^{\nu -1}}}} \\  \\  \\
\displaystyle =g_{j,0}(x,x-y)+\ M_j(x)\ \
\dfrac{\Omega _{\nu -1}}{{{(2\pi )^{\nu -1}}}}(\ln \vert x-y\vert
^{-1}+{\cal K%
}_{\nu -1}) \\  \\  \\
\displaystyle +\alpha _j^{\nu -j-2}(x,-1)+O(\vert x-y\vert ). \\
\end{array}
\end{equation}
The terms \\ \\
\begin{equation}
\frac{\partial _t^nA(x,0)\ }{n!}\int \left( \int_T^\infty t^n\
D_{-1-j}(x,t;\xi ,t;-1)\ dt\right) \ e^{i(x-y)\xi }{\frac{{d\xi
}}{{(2\pi
)^{\nu -1}}}}
\end{equation}
in (\ref{14}) have a limit when $y\rightarrow x$. Then, from (%
\ref{13}) to (\ref{16}) we obtain
\begin{equation}
\label{17}
\begin{array}{c}
\displaystyle \int_0^TA(x,t)\int D_{-1-j}(x,t;\xi ,t;-1)\
e^{i(x-y)\xi }\
{\dfrac{{d\xi }}{{(2\pi )^{\nu -1}}}}\ dt \\  \\  \\
\displaystyle =\sum\limits_{n=0}^{\nu -j-3}
\dfrac{\partial _t^nA(x,0)\ }{n!}g_{j,j+n-\nu +2}(x,x-y)\\ \\  \\
\displaystyle +\sum\limits_{n=0}^{\nu -j-3}
\dfrac{\partial _t^nA(x,0)\ }{n!}\left[ \alpha _j^n(x,-1)-\dfrac
1{-j-n+\nu
-2}\ \beta _j^n(x,-1)\right.  \\  \\  \\
\displaystyle -\left. \int \int_T^\infty t^n\ D_{-1-j}(x,t;\xi
,t;-1)\ dt\
{\dfrac{{d\xi }}{{(2\pi )^{\nu -1}}}}\right]  \\  \\  \\
\displaystyle +\dfrac{\partial _t^{\nu -j-2}A(x,0)}{(\nu
-j-2)!}\left[ g_{j,0}(x,x-j)+\ \
M_j(x)\ \dfrac{\Omega _{\nu -1}}{{{(2\pi )^{\nu -1}}}}(\ln \vert
x-y\vert
^{-1}+{\cal K}_{\nu -1})\right]  \\  \\  \\
\displaystyle +\dfrac{\partial _t^{\nu -j-2}A(x,0)}{(\nu
-j-2)!}\left[ \alpha _j^{\nu
-j-2}(x,-1)-\int \int_T^\infty t^{\nu -j-2}\ D_{-\nu +1}(x,t;\xi
,t;-1)\ dt\
{\dfrac{{d\xi }}{{(2\pi )^{\nu -1}}}}\right]  \\  \\  \\
\displaystyle +\int \left( \int_0^T\epsilon _j(x,t)\ t^{\nu -1-j}\
D_{-1-j}(x,t;\xi
,t;-1)\ dt\right) \ {\dfrac{{d\xi }}{{(2\pi )^{\nu -1}}}}+O(x-y),
\end{array}
\end{equation}
and then, comparing expressions (\ref{12}) and (\ref{17}), Lemma 4
follows. \qed
\bigskip\

Finally, in order to get part b) of Theorem 1 we write the equality
in Lemma 3 as
\begin{equation}
\lim _{z\rightarrow -1}\ R(x,t;z)=\lim _{y\rightarrow x}\ S(x,y,t)
\end{equation}
and evaluate the integrals $\int_0^TA(x,t)\ R(x,t;z)\ dt$ and $%
\int_0^TA(x,t)\ S(x,y,t)\ dt.$
 \eject
For the first one, we have  \\  \\
\begin{equation}
\label{18}
\begin{array}{c}
\displaystyle \int_0^TA(x,t)\left[ J_z(x,t;x,t)+\frac
1{z+1}\int_{\vert (\xi ,\tau )\vert
=1}c_{-\nu }(x,t;\xi ,\tau ;0)
\dfrac{\ d\sigma _{\xi ,\tau }}{(2\pi )^\nu }\right.  \\  \\  \\ \\
\displaystyle +\left. \int_{\vert (\xi ,\tau )\vert =1}\frac i{2\pi
}\int
\dfrac{\ln \lambda }\lambda c_{-\nu }(x,t;\xi ,\tau ;\lambda )\
d\lambda \
\dfrac{\ d\sigma _{\xi ,\tau }}{(2\pi )^\nu }\right] \ dt \\  \\  \\
\\
\displaystyle =-\sum\limits_{j=0}^{\nu -1}\int_0^TA(x,t)\ \int
D_{-1-j}(x,t;\xi ,t;z)\
{\dfrac{{d\xi }}{{(2\pi )^{\nu -1}}}}\ dt\\  \\  \\  \\
\displaystyle +\int_0^TA(x,t)\ R(x,t;z)\ dt \\   \\  \\
\\
\displaystyle =\sum\limits_{j=0}^{\nu -2}\frac 1{z+1}
\dfrac{\partial _t^{\nu -j-2}A(x,0)}{(\nu -j-2)!}\int_{\vert \xi
\vert
=1}\int_0^\infty t^{\nu -j-2}D_{-1-j}(x,t;\xi ,t;-1)\ \ dt\ \dfrac{\
d\sigma
_\xi }{(2\pi )^{\nu -1}} \\  \\  \\ \\
\displaystyle +\sum\limits_{j=0}^{\nu -2}
\dfrac{\partial _t^{\nu -j-2}A(x,0)}{(\nu -j-2)!}\int_{\vert \xi
\vert
=1}\int_0^\infty t^{\nu -j-2}\partial _zD_{-1-j}(x,t;\xi ,t;-1)\ \
dt\
\dfrac{\ d\sigma _\xi }{(2\pi )^{\nu -1}} \\  \\  \\ \\
\displaystyle -\sum\limits_{j=0}^{\nu -1}\psi _j(x,z)+\int_0^TA(x,t)\
\ R(x,t;z)\ dt.
\end{array}
\end{equation}
For the integral involving $\ S(x,y,t),$ we have
\begin{equation}
\label{19}
\begin{array}{c}
\displaystyle \int_0^TA(x,t)\left[ G_B(x,t;y,t)\right.  \\
\\
\left. -\sum\limits_{l=-(\nu -1)}^0h_l(x,t;x-y;0)-\ M(x,t)
\dfrac{\Omega _\nu }{(2\pi )^\nu }\left( \ln \vert x-y\vert
^{-1}+{\cal K}_\nu
\right) \ \right] dt \\  \\
\displaystyle =-\sum\limits_{j=0}^{\nu -1}\int_0^TA(x,t)\ \int
D_{-1-j}(x,t;\xi ,t;-1)\
e^{i(x-y)\xi }\
{\dfrac{{d\xi }}{{(2\pi )^{\nu -1}}}}\ dt \\  \\
\displaystyle +\int_0^TA(x,t)\ \ S(x,y,t)\ dt \\
\\
\displaystyle =-\sum\limits_{j=0}^{\nu -2}\sum\limits_{n=0}^{\nu
-j-2}
\dfrac{\partial _t^nA(x,0)}{n!}g_{j,j+n+\nu -2}(x,x-y) \\  \\
\displaystyle -\sum\limits_{j=0}^{\nu -2}
\dfrac{\partial _t^{\nu -j-2}A(x,0)}{(\nu -j-2)!}\ \
M_j(x)\dfrac{\Omega
_{\nu -1}}{(2\pi )^{\nu -1}}\left( \ln \vert x-y\vert ^{-1}+{\cal
K}_{\nu
-1}\right)  \\  \\
\displaystyle -\sum\limits_{j=0}^{\nu -1}\varphi
_j(x,y)+\int_0^TA(x,t)\ \ S(x,y,t)\ dt.
\end{array}
\end{equation}
Then, taking into account that the last terms in (\ref{18}) and
(\ref{19}) satisfy
 \begin{equation}
\begin{array}{c}
\displaystyle \lim \limits_{z\rightarrow -1}\left(
-\sum\limits_{j=0}^{\nu -1}\psi
_j(x,z)+\int_0^TA(x,t)\ \ R(x,t;z)\ dt\right) \\
\\
\displaystyle =\lim \limits_{y\rightarrow x}\left(
-\sum\limits_{j=0}^{\nu -1}\varphi
_j(x,y)+\int_0^TA(x,t)\ \ S(x,y,t)\ dt\right) ,
\end{array}
\end{equation}
we obtain part b) of Theorem 1. Notice that, for $\vert \xi \vert
\geq 1,$%
\begin{equation}
D_{-1-j}(x,t;\xi ,t;-1)=\ \ \ \tilde d_{-1-j}(x,t;\xi ,t;0)\
\end{equation}
and
\begin{equation}
\partial _zD_{-1-j}(x,t;\xi ,t;-1)=\frac i{2\pi }\int_\Gamma
\dfrac{\ln
\lambda }\lambda \ \ \ \tilde d_{-1-j}(x,t;\xi ,t;\lambda )\ d\lambda
{}.
\end{equation}

The proof of c) is similar to the one of b), and even simpler because
in
this case the parametrix in (\ref{100}) does not include terms of the
form
Op'($\theta ,$d$_{-1-j}).$ \qed

\bigskip
Awful as it looks, (\ref{rptmd}) is not so complicated: In the
first
place, all terms can be systematically evaluated. Moreover, the terms
containing $h_l$ subtract the singular part of the Green function in
the
interior of the manifold (see (\ref{?})) and can, thus, be easily
identified
from the knowledge of $G_B$. $R(x,t,y,t)$, the regular part so
obtained, is
still nonintegrable near the boundary. Those terms containing
$g_{j,l}$
subtract the singular part of the integrals $\int_0^T\ t^n\
R(x,t,y,t)\ dt$
(see (\ref{Asn})). Finally, the terms containing $c_{-\nu }$ and
$d_{-\nu
+1} $ arise as a consequence of having replaced an analytic
regularization
by a {\it point splitting } one.

Even though Seeley's coefficients $c$ and $\tilde d$ are to be
obtained
through an iterative procedure, which can make their evaluation a
tedious
task, in the cases of physical interest only the few
first of them are needed.
In fact, for the two dimensional example in the next section we will
only need two such coefficients.


\section{Two dimensional Dirac Operator on a disk.}

In this section, we will use the method previously discussed to
evaluate the
determinant of the operator $D=\ \not \!\!\!i\partial +\not \!\!A$
acting on functions
defined on a two dimensional disk of radius $R.$ A family of local
bag-like \cite{Bag}
elliptic boundary conditions will be assumed.

We take $A_\mu $ to be an Abelian field in the Lorentz gauge; as it
is well
known, it can be written as $A_\mu =\epsilon _{\mu \nu }\ \partial
_\nu \phi
\ (\epsilon _{01}=-\epsilon _{10}=1)$. For $\phi $ we choose a smooth
bounded
function $\phi =\phi (r)$ .
Notice that, with these assumptions, $A_r=0$ and $A_\theta
(r)=-\partial
_r\phi (r).$

We call
\begin{equation}
\label{3.2}\Phi =\oint_{r=R}\ A_\theta \ R\ d\theta =-2\pi R\
\partial
_r\phi (r) \vert _{r=R}.
\end{equation}

Our convention for two dimensional Dirac matrices is as in
(\ref{gm}):

\begin{equation}
\label{gammas}\gamma _0=\left(
\begin{array}{cc}
0 & 1 \\
1 & 0
\end{array}
\right) ~,\qquad \gamma _1=\left(
\begin{array}{cc}
0 & -i \\
i & 0
\end{array}
\right) ~,\qquad \gamma _5=\left(
\begin{array}{cc}
1 & 0 \\
0 & -1
\end{array}
\right) ,
\end{equation}
which satisfy:
\begin{equation}
\gamma _\mu \gamma _\nu \ =\delta _{\mu \nu }\ I\ +\ i\ \epsilon
_{\mu
\nu }\
\gamma _5.
\end{equation}
Therefore, the free Dirac operator can be written as:
\begin{equation}
i\not \! \partial =i\ (\gamma _0\ \partial _0+\gamma _1\ \partial
_1)=2i\
\left(
\begin{array}{cc}
0 & \frac \partial {\partial X} \\
\frac \partial {\partial X^{*}} & 0
\end{array}
\right) ,
\end{equation}
where $X=x_0+i\ x_1$ and $X^{*}=x_0-i\ x_1.$

Or, in polar coordinates:
\begin{equation}
i\not \! \partial =i(\gamma _r\ \partial _r+\frac 1r\gamma _\theta \
\partial
_\theta ),
\end{equation}

with
\begin{equation}
\label{gampol}\gamma _r=\left(
\begin{array}{cc}
0 & e^{-i\theta } \\
e^{i\theta } & 0
\end{array}
\right) ,\qquad \gamma _\theta =\left(
\begin{array}{cc}
0 & -ie^{-i\theta } \\
ie^{i\theta } & 0
\end{array}
\right) .
\end{equation}

With these conventions, the full Dirac operator can be written as:
\begin{equation}
\label{op}D=e^{-\gamma _5\phi (r)\ }i\not \! \partial \ e^{-\gamma
_5\phi (r)}.
\end{equation}

\bigskip

Now, in order to perform our calculations, we consider the family of
operators:
\begin{equation}
\label{opalf}D_\alpha =i\not \! \partial +\alpha \not \! \!
A=e^{-\alpha \gamma
_5\phi (r)\;\ }i\not \! \partial \ e^{-\alpha \gamma _5\phi (r)},\
\rm{with }\
0\leq \alpha \leq 1,
\end{equation}
which will allow us to go smoothly from the free to the full Dirac
operator.
If we call
\begin{equation}
W(\alpha )=\ln \ Det(D_\alpha )_B,
\end{equation}
where $B$ represents the  elliptic boundary
condition, we have
\begin{equation}
\dfrac \partial {\partial \alpha }W(\alpha )=\begin{array}{c}
\\
F.P. \\
^{_{z=0}}
\end{array}
\left[ Tr\left( \not\!\! A\
(D_\alpha )_B^{-z-1}\right) \right].
\end{equation}
{}From the Theorem in the previous section we get:
\begin{equation}
\label{ptmd}
\begin{array}{c}
\displaystyle \dfrac \partial {\partial \alpha }W(\alpha )=\frac
1{(2\pi )^2}\ tr\left\{
\int \ \lim \limits_{y\rightarrow x}\left[ \int \left[
\not \! \! A(t)\left( 4\pi ^2G_{B}(x,t,y,t)\right. \right. \right.
\right.
\\  \\
\displaystyle -\frac 1{\vert x-y\vert }\int\  e^{i\xi
\frac{(x-y)}{\vert x-y\vert }}\ c_{-1}(x,t;\frac{(\xi ,\tau )}{\vert
(\xi ,\tau
)\vert };0) \ d\xi \ d\tau \ \\  \\
\displaystyle -\int_{\vert (\xi ,\tau )\vert \geq 1}\ e^{i\xi (x-y)}\
c_{-2}(x,t;\xi ,\tau ;0)\ d\xi \ d\tau
\\   \\
\displaystyle \biggl. -\int\frac i{2\pi }\int_\Gamma \
\frac{\ln \lambda}\lambda \  c_{-2}(x,t;\frac{(\xi ,\tau )}{%
\vert (\xi ,\tau )\vert };\lambda )\ d\lambda\
d\sigma_{\xi,\tau}\biggr)\\  \\
\displaystyle +
\not \! \! A(0)\ \biggl( \int_{\vert \xi \vert \geq 1}\ e^{i\xi
(x-y)}\ \tilde
d_{-1}(x,t;\xi ,t;0)\ d\xi\biggr. \\  \\
\displaystyle \left. \left. \biggl. \biggl. +\ \int \frac i{2\pi
}\int_\Gamma \frac{\ln \lambda}\lambda \ \tilde d_{-1}(x,t;
\frac{(\xi )}{\vert \xi \vert },t;\lambda )\ d\lambda\ d\sigma_\xi
\biggr) \biggr]dt \right]dx \right\} ,
\end{array}
\end{equation}
where the Fourier transforms of $c_{-2}$ and $\tilde d_{-1}$ have
been left
explicitly indicated, instead of using the results of Lemma 1.

Now, the coefficients $c$ and $\tilde d$ in the previous equation are
those
appearing in the asymptotic expansion of the resolvent  $(D_\alpha
-\lambda I$
$)^{-1}$.

{}From (\ref{op}), the symbol of $(D_\alpha -\lambda I)$ is:
\begin{equation}
\label{sim}
\begin{array}{c}
\sigma (D_\alpha -\lambda I)=(-
\not  \xi -\lambda I)+\alpha \not \! \! A \\  \\
\displaystyle =a_1(\theta ,t,\xi ,\tau ,\lambda )+a_0(\theta ,t,\xi
,\tau ,\lambda ),
\end{array}
\end{equation}
where
\begin{equation}
\begin{array}{c}
a_1=-
\not  \xi -\lambda I, \\  \\
a_0=\alpha \not \! \! A.
\end{array}
\end{equation}
The c -coefficients can then be obtained from (\ref{7}) which,
written in
detail, gives
\begin{equation}
\begin{array}{c}
\displaystyle c_{-1}=a_1^{-1} \\
\\
\displaystyle a_1c_{-1-j}+\sum \left[ \left( D_{\xi ,\tau }\right)
^\beta a_k\right]
\left( iD_{x,t}\right) ^\beta c_{-1-l}/\beta !=0,\ \ \ j=1,2,...
\end{array}
\end{equation}
where the sum is taken over all $l<j$ and $k-\vert \beta \vert
-1-l=-j.$

So, the required Seeley's c-coefficients are given by \cite{Annals} :
\begin{equation}
\label{c}
\begin{array}{c}
\displaystyle c_{-1}=\frac 1{(\lambda ^2-\xi ^2-\tau ^2)}(
\not \! \xi -\lambda I), \\  \\
\displaystyle c_{-2}=\frac \alpha {(\lambda ^2-\xi ^2)^2}(2\lambda
\xi _\mu A_\mu
I-(\lambda ^2-\xi ^2)\not \! \! A-2\xi _\mu A_\mu \not \! \xi ),
\end{array}
\end{equation}
where $\not \xi =\xi \gamma _\theta +\tau \gamma _t$ .

As regards the boundary contributors to the parametrix, i.e., the
coefficients $d_{-1-j}$,
they are the solutions of (\ref{9}-\ref{12}).

Now, from (\ref{11.1})

\begin{equation}
\begin{array}{c}
\displaystyle a^{(0)}=a_0\vert _{t=0}=\alpha
\not \! \! A\vert _{r=R}, \\  \\
\displaystyle a^{(1)}=a_1\vert _{t=0}=-\lambda I-\xi \gamma _\theta
+i\gamma _t\dfrac
\partial {\partial t}.
\end{array}
\end{equation}
In our case, the equation to be solved is
\begin{equation}
\label{recast}a^{(1)}d_{-1}=(-\lambda I-\xi \gamma _\theta +i\gamma
_t\partial _t)d_{-1}=0,
\end{equation}
with boundary conditions
\begin{equation}
\label{bcd1}b_0\ d_{-1}=b_0\ c_{-1}\ \ \rm{\ at }\ t=0,
\end{equation}
plus the vanishing of $\ d_{-1}$ as $t\rightarrow +\infty
.\;$(\ref{recast})
can be recast in the form
\begin{equation}
\label{eqd1}\partial _td_{-1}=-M\ d_{-1,}
\end{equation}
where $M=\xi \gamma _5+i\lambda \gamma _t$ . It can be easily
verified that
\begin{equation}
\begin{array}{c}
tr(M)=0, \\
\\
M^2=(\xi ^2-\lambda ^2)I.
\end{array}
\end{equation}
So, $M$ has eigenvalues $\pm \sqrt{\xi ^2-\lambda ^2}$, corresponding
to the
eigenvectors
\begin{equation}
u_{\pm }=\left(
\begin{array}{c}
ie^{-i\theta }(\xi \pm
\sqrt{\xi ^2-\lambda ^2}) \\ \lambda
\end{array}
\right) .
\end{equation}
Now, the general solution to (\ref{eqd1}) is:
\begin{equation}
d_{-1}(x,t;\xi ,\tau ;\lambda )=e^{-tM}\ d_{-1}(x,0;\xi ,\tau
;\lambda ).
\end{equation}
Since $d_{-1}\rightarrow 0$ for $t\rightarrow \infty $ , we get
\begin{equation}
\label{d}d_{-1}(x,t;\xi ,\tau ;\lambda )=e^{-t\sqrt{\xi ^2-\lambda
^2}%
}u_{+}\otimes \left( \QATOP{f}{g}\right) ^{\dagger },
\end{equation}
where the vector $\left( \QATOP{f}{g}\right) $ must be determined
from the
boundary condition at $t=0~(r=R),$ given by (\ref{bcd1}). After this
brief review of the general points, let us now go to the specific
calculations.

We now consider a parametric family of bag-like local
boundary conditions leading to an elliptic boundary problem.
According to the discussion leading to (\ref{betas}), we choose the
matrix
\begin{equation}
\label{b0}
b_0=\left( 1,w\ e^{-i\theta }\right) ,
\end{equation}
with $w$ a nonzero complex constant.

We define the operator $(D_\alpha )_B$ as the differential operator
in (%
\ref{op}) , acting on the dense subspace of functions satisfying:
\begin{equation}
\label{lbc}
B\ \psi \equiv b_0 \psi \vert _{t=0}=0.
\end{equation}

Notice that these boundary conditions reduce to those of an MIT bag
\nobreak \cite{Bag} when $
w=\pm 1$ .

\subsection{Zero modes}

We will here show that, with these boundary conditions, the operator
is
invertible. From (\ref{op}), we get:

\begin{equation}
D_\alpha \ \psi =0\Rightarrow \not \! \partial \ e^{-\alpha \gamma
_5\phi (r)}\
\psi =0,
\end{equation}
or, equivalently:
\begin{equation}
\left(
\begin{array}{cc}
0 & e^{-i\theta }(\partial _r-\frac ir\partial _\theta ) \\
e^{i\theta }(\partial _r+\frac ir\partial _\theta ) & 0
\end{array}
\right) \left(
\begin{array}{cc}
e^{-\alpha \phi (r)} & 0 \\
0 & e^{\alpha \phi (r)}
\end{array}
\right) \left( \QATOP{\varphi (r,\theta )}{\chi (r,\theta )}\right)
=0.
\end{equation}
Now, we introduce the expansions:

\begin{equation}
\begin{array}{c}
\displaystyle \varphi (r,\theta )=\sum\limits_{n=-\infty }^\infty
\varphi _n(r)\
e^{in\theta }, \\
\\
\displaystyle \chi (r,\theta )=\sum\limits_{n=-\infty }^\infty \chi
_n(r)\ e^{in\theta }.
\end{array}
\end{equation}
The solutions are thus given by:
\begin{equation}
\begin{array}{c}
\displaystyle \varphi _n(r)=a_{n\ }\ r^n\ e^{\alpha \phi (r)}, \\
\\
\displaystyle \chi _n(r)=b_{n\ }\ r^{-n}\ e^{-\alpha \phi (r)},
\end{array}
\end{equation}
where the coefficients $a_n$ and $b_n$ are to be determined from the
boundary conditions at $r=0$ and $r=R$ . The requirement of
normalizability
 implies that :
\begin{equation}
\begin{array}{c}
a_n=0,\ \
\rm{ for }\ n<0, \\ b_n=0,\ \ \rm{ for }\ n>0.
\end{array}
\end{equation}
At $r=R$, we get, from (\ref{lbc}):
\begin{equation}
b_{n+1}=-\frac{a_n}w\ R^{2n+1}\ e^{2\alpha \phi (R)},
\end{equation}
which requires

\begin{equation}
\begin{array}{c}
a_n=0,\ \
\rm{ for }\ n\geq 0, \\ b_n=0,\ \ \rm{ for }\ n\leq 0.
\end{array}
\end{equation}
So, for these local boundary conditions, there are no normalizable
zero
modes. Notice that these are not the most general local elliptic
boundary
conditions. In fact, zero modes would in general arise if one allowed
$w$ to
depend on $\theta $.

\subsection{ Computation of $d_{-1}$}

Restricting ourselves to the case of $\theta $-independent parameter
$w$,  we
now go back to (\ref{d}), from which:
\begin{equation}
d_{-1}\vert _{r=R}=u_{+}\otimes \left( \QATOP{f}{g}\right) ^{\dagger
},
\end{equation}
and determine $\left( \QATOP{f}{g}\right) ^{\dagger }$ from the
boundary
condition:
\begin{equation}
b_0\ d_{-1}=b_0\ c_{-1},\ \ \rm{ at \ }r=R,
\end{equation}
which, written in an explicit manner reads:

\begin{equation}
\{(1,w\ e^{-i\theta })\ .\ u_{+}\}\left( \QATOP{f}{g}\right)
^{\dagger
}=(1,w\ e^{-i\theta })\ c_{-1}.
\end{equation}

{}From the expression for $c_{-1}$ given in  (\ref{c}), it turns out
that:
\begin{equation}
\begin{array}{c}
\left( \QDATOP{f}{g}\right) ^{\dagger }=\dfrac{e^{i\theta }}{(\xi
^2+\tau
^2-\lambda ^2) (\lambda w+i\xi +i\sqrt{\xi ^2-\lambda ^2)}}\hfill
\\ \\ \qquad \qquad \qquad \qquad \times  \left(\ \lambda +w\
(-i\xi +\tau ) \qquad e^{-i\theta }(i\xi +\tau +\lambda \ w)\
\right).
\end{array}
\end{equation}
Replacing this expression into (\ref{d}), we obtain:
\begin{equation}
\begin{array}{c}
d_{-1}=
\dfrac{e^{-t\sqrt{\xi ^2-\lambda ^2}}}{(\xi ^2+\tau ^2-\lambda ^2)\
(w\lambda
+i\xi +i\sqrt{\xi ^2-\lambda ^2})}\hfill
 \\  \\  \qquad \qquad \times \left(
\begin{array}{cc}
\scriptstyle i(\xi +\sqrt{\xi ^2-\lambda ^2})\ (\lambda -w(i\xi -\tau
)) & \ \ \ \scriptstyle ie^{-i\theta
}\ (\xi +
\sqrt{\xi ^2-\lambda ^2})\ (w\lambda +i\xi +\tau ) \\  \\
\scriptstyle \ \lambda \ e^{i\theta
}\ (\lambda -w(i\xi -\tau )) & \scriptstyle \ \lambda \ (w\lambda
+i\xi +\tau )
\end{array}
\right) .
\end{array}
\end{equation}

Now, taking into account (\ref{`d}), i.e.,
 \begin{equation}
\tilde d_{-1}=-\oint_{\Gamma ^{-}}d\tau \ e^{-i\tau u}d_{-1}(\theta
,t,\xi
,\tau ,\lambda ),
\end{equation}
we finally get:
\begin{equation}
\label{dm1til}
\begin{array}{c}
\tilde d_{-1}=\pi i
\dfrac{e^{-(u+t)\sqrt{\xi ^2-\lambda ^2}}}{\sqrt{\xi ^2-\lambda
^2}(iw\lambda
-\xi -\sqrt{\xi ^2-\lambda ^2})}\ \ \hfill \\  \\
\times \left(
\begin{array}{cc}
\scriptstyle (\xi +\sqrt{\xi ^2-\lambda ^2})\ (i\lambda +w(\xi
+\sqrt{\xi ^2-\lambda ^2}
)) \ \ & \scriptstyle e^{-i\theta }(\xi +
\sqrt{\xi ^2-\lambda ^2})\ (iw\lambda -\xi +\sqrt{\xi ^2-\lambda ^2})
\\  \\
\scriptstyle -i\lambda \ e^{i\theta }(i\lambda -w(\xi +\sqrt{\xi
^2-\lambda ^2})) &
\scriptstyle -i\lambda \ (iw\lambda -\xi +\sqrt{\xi ^2-\lambda ^2})
\end{array}
\right) .
\end{array}
\end{equation}

\subsection{Calculation of the Green function}

We are looking for the function $G_{B}(x,y)$ satisfying:
\begin{equation}
\begin{array}{c}
D_\alpha \ G_{B}(x,y)=\delta (x,y), \\
\\
B\ G_{B}(x,y)\vert _{x\in \partial \Omega }=0,
\end{array}
\end{equation}
where $D_\alpha $ and $B,$ are given by equations (\ref{op}) and
(\ref{lbc})
respectively. Now, some notation is in order:
\begin{equation}
\begin{array}{c}
x=(x_0 , x_1)=(r\cos \theta ,r\sin \theta ), \\
X=x_0+i\ x_1=r\ e^{i\theta }, \\
\\
y=(y_0 , y_1)=(\rho \cos \varphi ,\rho \sin \varphi ),\  \\
Y=y_0+i\ y_1=\rho \ e^{i\varphi }.
\end{array}
\end{equation}
It is easy to see that $\ G_{B}(x,y)$ can be written as
\begin{equation}
\ G_{B}(x,y)=e^{\alpha \gamma _5\phi (r)}\left[
G_0(x,y)+h(x,y)\right]
e^{\alpha \gamma _5\phi (\rho )},
\end{equation}
where
\begin{equation}
\ G_0(x,y)=\frac 1{2\pi i}\frac{\not\! x-\not\! y}{(x-y)^2}
\end{equation}
 is the Green function of the operator $\not \!\!\!\!i\partial$ and
$h(x,y)$ is a solution of the homogeneous equation
\begin{equation}
i\not \! \partial \ h(x,y)=0,
\end{equation}
to be determined through the boundary condition (\ref{b0}). Due to
the
geometry of the problem, one can make the ansatz
\begin{equation}
h(x,y)=\ G_0(x,\tilde y)\ H(y),
\end{equation}
where $\tilde y$ is related to $y$ trough the inversion $\tilde
y=y\frac{R^2}{\rho ^2}$ . Then, taking into account
that
\begin{equation}
\frac 1r\ \gamma _r\ G_0(\tilde x,y)=-G_0(x,\tilde y)\ \frac 1\rho \
\gamma
_\rho ,
\end{equation}
one finds:
\begin{equation}
H(y)=e^{2\alpha \gamma _5\phi (r)}\ \frac R\rho \ \frac{\left[
(1+w^2)\
I+(1-w^2)\ \gamma _5\right] }{2w}\ \gamma _\rho .
\end{equation}
Thus, the relevant Green function is given by
\begin{equation}
\label{green1}G_{B}(x,y)=\frac 1{2\pi i}\left(
\begin{array}{cc}
\frac{R\ w\ e^{\alpha (\phi (x)+\phi (y)-2\phi (R))}}{XY^{*}-R^2} &
\frac{%
e^{\alpha (\phi (x)-\phi (y))}}{X-Y} \\  \\ \frac{e^{-\alpha (\phi
(x)-\phi
(y))} }{(X-Y)^{*}} & \frac{R\ e^{-\alpha (\phi (x)+\phi (y)-2\phi
(R))}}{w\
(XY^{*}-R^2)^{*}}
\end{array}
\right) .
\end{equation}

\subsection{Evaluation of the determinant}

With these elements at hand, we now go back to (\ref{ptmd}), and
perform the calculation of the determinant .

{}From (\ref{green1}), one can see that
\begin{equation}
G_{B}(\theta ,r,\varphi ,r)
\begin{array}{c}
\\
\sim \\
\scriptstyle \varphi \rightarrow \theta
\end{array}
\rm{diagonal\ matrix }+\frac 1{2\pi i\ r\ (\theta -\varphi )}\ \gamma
_\theta .
\end{equation}

When replaced into (\ref{ptmd}), we get for the first term in
the r.h.s.

\begin{equation}
tr\left\{ A_\theta \ \gamma _\theta \ G_B(\theta ,r,\varphi
,r)\right\}
\begin{array}{c}
\\
\sim \\
\scriptstyle \varphi \rightarrow \theta
\end{array}
\frac{A_\theta }{\pi i\ r\ (\theta -\varphi )}.
\end{equation}

For the second term in  (\ref{ptmd})
\begin{equation}
-\frac 1{4\pi ^2\vert x-y\vert }\int d\xi \ d\tau \ e^{i\xi
\frac{(x-y)}{\vert
x-y\vert }}\ c_{-1}(x,t;\frac{(\xi ,\tau )}{\vert (\xi ,\tau )\vert
};\lambda
=0)
\begin{array}{c}
\\
\sim \\
\scriptstyle \varphi \rightarrow \theta
\end{array}
\frac{-1}{2\pi i\ r\ (\theta -\varphi )}\ \gamma _\theta ,
\end{equation}
which exactly cancels the singularity of the Green function. (Notice
that
this Fourier transform must be understood in the sense of
distributions). Therefore, the
contribution of the first two terms in  (\ref{ptmd}) vanishes.

As regards the third term,
 \begin{equation}
\label{w3}
\begin{array}{c}
\displaystyle \frac{-1}{(2\pi )^2}\ tr\int \lim \limits_{y\rightarrow
x}\ \not\!\!
A(t)\int_{\vert (\xi ,\tau )\vert \geq 1}e^{i\xi (x-y)}\
c_{-2}(x,t;\xi ,\tau ;0)\ d\xi \ d\tau\  dx\ dt\\  \\
\displaystyle=  \frac{-\alpha }{2\pi ^2}\ \lim \limits_{y\rightarrow
x}\ \int  A_\theta
^2\ d^2x\ \int_{\vert (\xi ,\tau )\vert \geq 1} e^{i\xi (x-y)\
}\frac{%
(\tau ^2-\xi ^2)}{(\xi ^2+\tau ^2)^2}\ d\xi \ d\tau \\  \\
\displaystyle=
\frac{-\alpha }\pi \int A_\theta ^2\ d^2x\ \ \lim
\limits_{y\rightarrow
x}\int_{\vert x-y\vert }^\infty J_2(u) \ \  \frac{du}u\ =\
\frac{-\alpha }{2\pi }\int A_\nu \ A_\nu\ d^2x\  .
\end{array}
\end{equation}
where $J_2(u)$ is the Bessel function of order two.

Now, the fourth term in ( \ref{ptmd}) is:
\begin{equation}
\label{w4}
\begin{array}{c}
\displaystyle\frac{-1}{(2\pi )^2}\ tr\int \not \! \! A(t)\int \frac
i{2\pi }\int_\Gamma \ \ln \lambda \ \ c_{-2}(x,t;
\frac{(\xi ,\tau )}{\vert (\xi ,\tau )\vert};\lambda ) \ \
\frac{d\lambda }
\lambda \ d\sigma_{\xi,\tau}\ dx\ dt\\\  \\
\displaystyle=
\frac{-i\alpha }{4\pi ^3}\int A_\theta
^2\ d^2x\ \int_\Gamma \ \ \frac{\ln \lambda }{(\lambda ^2-1)^2}
\int \ (1-\lambda ^2-2\xi ^2)\ d\sigma_{\xi,\tau}\  \frac{d\lambda }
\lambda \ \\  \\
\displaystyle=\frac{i\alpha }{2\pi ^2}\int A_\theta
^2\ d^2x\ \ 2\pi i\ \int_0^\infty \frac{%
\mu \ d\mu }{(\mu ^2+1)^2} \ =\ \frac{-\alpha }{2\pi }\int A_\nu \
A_\nu\ d^2x\  .
\end{array}
\end{equation}
This term gives rise to a contribution identical to that of
(\ref{w3}).

The last term in (\ref{ptmd}) is

\begin{equation}
\label{last}
\begin{array}{c}
\displaystyle\frac i{(2\pi )^3}\ tr\int
\not \! \! A(0)\ \sum\limits_{\xi =\pm 1}\int_\Gamma \ \ln
\lambda \ \ \tilde d_{-1}(x,t;\frac \xi {\vert \xi \vert },t;\lambda
)
\ \frac{d\lambda }\lambda\ dx\ dt\ \\  \\
\displaystyle=\frac{i\ \Phi }{(2\pi )^2}\int_\Gamma
\ \frac{u\ \ln \lambda }{(1+\ u^2\ \lambda ^2)}[\lambda \
\sqrt{1+u^2}-i
\sqrt{1-\lambda ^2}]\ \frac{d\lambda }{\sqrt{1-\lambda ^2}} ,
\end{array}
\end{equation}
where $u=(1-w^2)\ /\ 2w.$ We choose the curve $\Gamma $ as in Fig.
\ref
{figura1}.

\begin{figure}
\begin{center}
\begin{picture}(150,150)(-75,-50)
\put (0,-50){\line(0,1){150}}
\put (-75,0){\line(1,0){150}}
\put (5,100){\vector(0,-1){60}}
\put (5,40){\line(0,-1){30}}
\put (5,0){\oval(20,20)[r]}
\put (5,-10){\line(-1,0){10}}
\put (-5,0){\oval(20,20)[l]}
\put (-5,10){\vector(0,1){60}}
\put (-5,70){\line(0,1){30}}
\put (20,20){\makebox(0,0){$\Gamma$}}

\end{picture}
\end{center}
\caption{ The contour $\Gamma$ }
\label{figura1}
\end{figure}

Therefore, (\ref{last}) reads
\begin{equation}
\label{w6}
\begin{array}{c}
\displaystyle -\dfrac \Phi {2\pi }\ u\ \int_0^\infty \dfrac 1{(1-\
u^2\ \mu
^2)}\left[ \mu \ \
\dfrac{\sqrt{\ 1\ +\ u^2}}{\sqrt{\ 1\ +\ \mu ^2}}\ -\ 1\right] d\mu \
\\  \\
=-\dfrac \Phi {4\pi }\left[ \ln \left(
\dfrac{\sqrt{\ 1\ +\ u^2}\ +\ u\ \sqrt{\ 1\ +\ \mu ^2}}{\sqrt{\ 1\ +\
u^2}\
-\ u\ \sqrt{\ 1\ +\ \mu ^2}}\ \ \dfrac{1\ -\ u\ \mu }{1\ +\ u\ \mu
}\right)
\right] _0^\infty \\  \\
=\dfrac{-\Phi }{4\pi }\ \ln \ w^2.
\end{array}
\end{equation}

Putting all pieces together ((\ref{w3}), (\ref{w4}) and (\ref{w6}),
we
finally find:
\begin{equation}
\begin{array}{c}
\displaystyle\ln Det (D)_B-\ln Det (\!i\not\!\partial )_B=-\frac
1{2\pi } \int_\Omega A_\nu \ A_\nu \ d^2x\ \
-\frac \Phi {4\pi }\ln \ w^2 \\  \\
\displaystyle=-\frac 1{2\pi }\int_\Omega  A_\nu \ A_\nu\  d^2x\ \
 -\ \frac 1{4\pi }\
\ln \ w^2\int_{\partial \Omega }A_{\nu }\ dx_\nu .
\end{array}
\end{equation}
The first term is the integral, restricted to the region $\Omega $,
of the
same density appearing in the well known case of the whole plane
\cite
{Schwinger}. The second term is well-defined for every $w\neq 0$, and
vanishes for a null total flux, $\Phi =0$. For $w=0$, $b_0$ in
(\ref{b0}) does not define an elliptic boundary problem, as discussed
in Section 2. It is also interesting to notice that this term
vanishes in the case of MIT bag boundary conditions, i.e., $w=\pm 1.$

This calculation is to be compared with the case of the compactified
plane \cite{Annals}, where the determinant can be expressed in terms
of just the kernel of the $z$-power of the operator analytically
extended to $z=0$, which is a local quantity. The presence of
boundaries makes the evaluation more involved, since even in simple
cases as the present (or the half plane treated in \cite{Jour.2}),
the knowledge of the Green function of the problem is needed.

\bigskip

{\it Acknowledgments}.  We are grateful to R.T. Seeley for useful
comments.
\bigskip

\bigskip
\eject


\begin{thebibliography}{99}
\bibitem[1]{APS}Atiyah, M.F., Patodi, V.K., Singer, I.M.: Math. Proc.
Camb. Phil. Soc. {\bf 77}, 43 (1975);  {\bf 78}, 405 (1975); {\bf
79}, 71 (1976).

\bibitem[2]{Booss-B}  B. Booss, B.,  Bleecker, D.D.: Topology and
Analysis. The Atiyah-Singer Index Formula and Gauge-$\!$Theoretic
Physics.
Springer-Verlag, New York, 1985.

\bibitem[3]{Booss-W}Booss, D., Wojciechowski, K.P. : Elliptic
Boundary Problems for Dirac Operators. Birkh\"auser, Boston, 1993.

\bibitem[4]{LibroAmarillo}
Calder\'on, A.P.: Lectures notes on
pseudodifferential operators and elliptic boundary value problems, I.
Publicaciones del I.A.M., Buenos Aires, 1976.

\bibitem[5]{Bag}Chodos, A., Jaffe, R.L., Johnson, K., Thorn, C.B.,
Weisskoft, V.F.: Phys. Rev. {\bf D9},
 3471 (1974); A. Chodos, R.L. Jaffe, K. Johnson and C.B. Thorn: Phys.
Rev. {\bf D10}, 2599 (1974).

\bibitem[6]{Jour.2}Falomir, H., Muschietti, M.A., Santangelo, E.M.:
J. Math. Phys. {\bf 31}, 989 (1990).

\bibitem[7]{Forman}Forman, R.: Functional Determinants and geometry.
Invent. math. {\bf 88}, 447 (1987).

\bibitem[8]{Annals}Gamboa Sarav\'\i , R.E., Muschietti, M.A.,
Schaposnik F.A., Solomin, J.E.:
Ann. Phys. {\bf 157}, 360 (1984).

\bibitem[9]{Jour.1} Gamboa Sarav\'\i , R.E., Muschietti, M.A.,
Schaposnik, F.A., Solomin, J.E.:
J. Math. Phys. {\bf 26}, 2045 (1985).

\bibitem[10]{Gil}Gilkey, P.B., Smith,L.: Communications on Pure and
Applied Mathematics, Vol. XXXVI, 85-131 (1983).

\eject

\bibitem[11]{SG}Grubb, G., Seeley, R.T.: Developpements Asymptotiques
pour l'operateur d'Atiyah-Patodi-Singer.
C.R. Acad. Sci. Paris, {\bf 317}, serie  I, p. 1123-1126 (1993);\\
Zeta and eta fuctions for Atiyah-Patodi-Singer operators, Copenh.
Univ., Math.  Dept, Prepint Ser., (1993), nro. 11;\\
Weakly parametric pseudodifferentials operators and
Atiyah-Patodi-Singer boundary ploblems, Copenh. Univ.,
Math. Dept, Prepint Ser., (1993), nro 15.

\bibitem[12]{Haw}Hawking, S.W.: Comm. Math. Phys. {\bf 55}, 133
(1977).

\bibitem[13]{Hor}H\"ormander, L.: The Analysis of Linear Partial
Differential Operators III, Pseudo-Differential Operators.
Springer-Verlag, Berlin Heidelberg, 1985.

\bibitem[14]{Moreno}Moreno, E.F.:  Phys. Rev. {\bf D48}, 921 (1993).

\bibitem[15]{NT}Ninomiya, M., Tan, C.I.:  Nucl. Phys. {\bf B257}, 199
(1985).

\bibitem[16]{Schwinger}Schwinger,J.: Phys. Rev. {\bf 82}, 664 (1951);
{\bf 128}, 2425 (1962).

\bibitem[17]{Seeley1}Seeley, R.T.: The resolvent of an elliptic
boundary problem. Am. J. Math. 91, 889-920 (1969).

\bibitem[18]{Seeley2}Seeley, R.T.: Analytic extension of trace
associated with elliptic boundary problems. Am. J. Math. 91, 963-983
(1969).

\bibitem[19]{Stein}Stein, E.M.: Singular Integrals and
Differentiability Properties of functions. Princeton University
Press,
Princeton, New Jersey, 1970.

\end{thebibliography}
\end{document}